\documentclass[11pt,twoside,leqno]{aomamlt2e} 
  \pageno{1305}
\received{June 11, 2003}
 
 \renewcommand{\theequation}{\thesection.\arabic{equation}}

\newtheorem{thm}{Theorem}[section]

\newtheorem{lemma}[thm]{Lemma}

\newtheorem{claim}[thm]{Claim}
\theorembodyfont{\upshape}
\newtheorem{rmk}[thm]{{\it Remark}}

\newcommand{\R}{\mathbb{R}}
\newcommand{\Z}{\mathbb{Z}}

\newcommand{\bdry}{\partial}
\renewcommand{\L}{\mathcal{L}}
\newcommand{\K}{\mathcal{K}}
\newcommand{\CC}{\mathcal{C}}

\newcommand{\tb}{tb}

\newcommand{\s}{\vskip.13in}
\newcommand{\n}{\noindent}

\begin{document}
\currannalsline{162}{2005} 

\title{Cabling and transverse simplicity}

\twoauthors{John B.\ Etnyre}{Ko Honda}

 \institution{University of Pennsylvania,
Philadelphia, PA  \\
\email{etnyre@math.upenn.edu}\\
\url{http://www.math.upenn.edu/$\sim$etnyre}\\
\\
\vglue-9pt
University of Southern California,
Los Angeles, CA\\
\email{khonda@math.usc.edu}\\
\url{http://rcf.usc.edu/$\sim$khonda}}%\char126

\centerline{\bf Abstract}
\vskip12pt

We study Legendrian knots in a cabled knot type.  Specifically, given
a topological knot type $\K$, we analyze the Legendrian knots in knot
types obtained from $\K$ by cabling, in terms of Legendrian knots in
the knot type $\K$.  As a corollary of this analysis, we show that the
$(2,3)$-cable of the $(2,3)$-torus knot is {\em not transversely
  simple} and moreover classify the transverse knots in this knot
type.  This is the first classification of transverse knots in a
non-transversely-simple knot type. We also classify Legendrian knots
in this knot type and exhibit the first example of a Legendrian knot
that does not destabilize, yet its Thurston-Bennequin invariant is not
maximal among Legendrian representatives in its knot type.

%%%%%%%%%%%%%%%%%%%%%%%%%%%%%%%%%%%%%%%%%%%%%%%%%%%%%%%%%%

%----------------------------------------------------------------------------
\section{Introduction}
%----------------------------------------------------------------------------

In this paper we continue the investigation of Legendrian knots in 
tight contact $3$-manifolds using $3$-dimensional contact-topological 
methods.  In \cite{EH1}, the authors introduced a general framework
for analyzing Legendrian knots in tight contact $3$-manifolds.  There we
streamlined the proof of the classification of Legendrian unknots,
originally proved by Eliashberg-Fraser in \cite{EF}, and gave a
complete classification of Legendrian torus knots and figure eight
knots.  In \cite{EH2}, we gave the first structure theorem for
Legendrian knots, namely the reduction of the analysis of connected
sums of Legendrian knots to that of the prime summands.  This yielded
a plethora of non-Legendrian-simple knot types. (A topological knot
type is {\em Legendrian simple} if Legendrian knots in this knot type
are determined by their Thurston-Bennequin invariant and rotation
number.) Moreover, we exhibited pairs of Legendrian knots in the same
topological knot type with the same Thurston-Bennequin and rotation
numbers, which required arbitrarily many stabilizations before they
became Legendrian isotopic (see \cite{EH2}).

The goal of the current paper is to extend the results obtained for
Legendrian torus knots to Legendrian representatives of cables of knot
types we already understand. On the way to this goal, we encounter the
{\em contact width}, a new knot invariant which is related to the
maximal Thurston-Bennequin invariant.  It turns out that the structure
theorems for cabled knots types are not as simple as one might expect,
and rely on properties associated to the contact width of a knot.
When these properties are not satisfied, a rather unexpected and
surprising phenomenon occurs for Legendrian cables. This phenomenon
allows us to show, for example, that the $(2,3)$-cable of the
$(2,3)$-torus knot is not transversely simple! (A topological knot
type is {\em transversely simple} if transverse knots in that knot
type are determined by their self-linking number.)  Knots which are
not transversely simple were also recently found in the work of Birman
and Menasco \cite{BM}.  Using braid-theoretic techniques they showed
that many three-braids are not transversely simple.  Our technique
should also provide infinite families of non-transversely-simple knots
(essentially certain cables of positive torus knots), but for
simplicity we content ourselves with the above-mentioned example.
Moreover, we give a complete classification of transverse (and
Legendrian) knots for the $(2,3)$-cable of the $(2,3)$-torus knot. 
This is the first classification of transverse knots in a
non-transversely-simple knot type.  

We assume that the reader has familiarity with \cite{EH1}.  In this
paper, the ambient $3$-manifold is the standard tight contact
$(S^3,\xi_{std})$, and all knots and knot types are {\em oriented}.
Let $\K$ be a topological knot type and $\L(\K)$ be the set of
Legendrian isotopy classes of $\K$.  For each $[L]\in \L(\K)$ (we
often write $L$ to mean $[L]$), there are two so-called {\em classical
  invariants}, the {\em Thurston-Bennequin invariant} $\tb(L)$ and the
{\em rotation number} $r(L)$.  To each $\K$ we may associate an
oriented knot invariant  
$$\overline{\tb}(\K)=\max_{L\in \L(\K)} \tb(L),$$ called the {\it 
maximal Thurston-Bennequin number}.   

A close cousin of $\overline{\tb}(\K)$ is another oriented knot 
invariant called the {\em contact width} $w(\K)$ (or simply the {\em 
width}) defined as follows:  First, an embedding $\phi: S^1\times 
D^2\hookrightarrow S^3$ is said to {\em represent} $\K$ if the core 
curve of $\phi(S^1\times D^2)$ is isotopic to $\K$.  (For notational 
convenience, we will suppress the distinction between $S^1\times D^2$ 
and its image under $\phi$.)  Next, in order to measure the {\em 
slope} of homotopically nontrivial curves on $\bdry (S^1\times D^2)$, 
we make a (somewhat nonstandard) oriented identification $\bdry 
(S^1\times D^2)\simeq \R^2/\Z^2$, where the meridian has slope $0$ and 
the longitude (well-defined since $\K$ is inside $S^3$) has slope 
$\infty$. We will call this coordinate system $\mathcal{C}_\K$. 
Finally we define 
$$w(\K)=\sup {1\over {\rm slope} (\Gamma_{\bdry (S^1\times D^2)})},$$
where the supremum is taken over  $S^1\times D^2\hookrightarrow S^3$ 
representing $\mathcal{K}$ with\break  $\bdry (S^1\times D^2) \mbox{ convex}.$

Note that there are several notions similar to $w(\K)$ --- see 
\cite{Colin}, \cite{Gay}.  The contact width clearly satisfies the following 
inequality: $$\overline{\tb}(\K)\leq w(\K) \leq 
\overline{\tb}(\K)+1.$$  In general, it requires significantly more 
effort to determine $w(\K)$ than it does to determine 
$\overline{\tb}(\K)$.   Observe that
$\overline{\tb}(\K)=-1$ and $w(\K)=0$ when $\K$ is the unknot.

\Subsec{Cablings and the uniform thickness property}
Recall that a {\em $(p,q)$-cable} $\K_{(p,q)}$ of a topological knot
type $\K$ is the isotopy class of a knot of slope $\frac{q}{p}$ on the
boundary of a solid torus $S^1\times D^2$ which represents $\K$, where
the slope is measured with respect to $\mathcal{C}_\K$, defined above.
In other words, a representative of $\K_{(p,q)}$ winds $p$ times
around the meridian of $\K$ and $q$ times around the longitude of
$\K$.  A {\em $(p,q)$-torus knot} is the $(p,q)$-cable of the unknot.    

One would like to classify Legendrian knots in a cabled knot type. This
turns out to be somewhat subtle and relies on the following key notion:

\demo{Uniform thickness property {\rm (UTP)}}  
Let $\K$ be a topological knot type.  Then $\K$ satisfies the {\it
  uniform thickness condition} or {\it is uniformly thick} if the
following hold: 

\begin{enumerate}
\item $\overline{\tb}(\K)=w(\K).$
\item Every embedded solid torus $S^1\times D^2\hookrightarrow S^3$
  representing $\K$ can be thickened to a {\em standard neighborhood}
  of a maximal $\tb$ Legendrian knot.  
\end{enumerate}

Here, a {\em standard neighborhood} $N(L)$ of a Legendrian knot $L$ is
an embedded solid torus with core curve $L$ and convex boundary $\bdry
N(L)$ so that $\#\Gamma_{\bdry N(L)}=2$ and $\tb(L)= {1\over
  {\rm slope} (\Gamma_{\bdry N(L)})}$.  Such a standard neighborhood
$N(L)$ is contact isotopic to any sufficiently small tubular
neighborhood $N$ of $L$ with $\bdry N$ convex and  
$\#\Gamma_{\bdry N}=2$. (See \cite{H1}.)  Note that, strictly speaking, 
Condition~2 implies Condition~1; it is useful to keep in mind,
however, that the verification of the UTP usually proceeds by
outlawing solid tori representing $\K$ with ${1\over
  {\rm slope} (\Gamma)}> \overline{\tb}(\K)$ and then showing that
solid tori with ${1\over {\rm slope} (\Gamma)} < \overline{tb}(\K)$
can be thickened properly.  We will often say that a solid torus $N$
(with convex boundary) representing $\K$ {\em does not admit a
  thickening}, if there is no thickening $N'\supset N$ whose
${\rm slope} (\Gamma_{\bdry N'})\not = {\rm slope} (\Gamma_{\bdry
  N})$.  

The reason for introducing the UTP is due (in part) to:

\begin{thm} \label{negative cables}
Let $\K$ be a knot type which is Legendrian simple and satisfies the
{\rm UTP}.  Then $\K_{(p,q)}$ is Legendrian simple and admits a
classification in terms of the classification of $\K$.  
\end{thm}

Of course this theorem is of no use if we cannot find knots satisfying
the UTP.  The search for such knot types has an inauspicious start as
we first observe that the unknot $\K$ {\em does not} satisfy the UTP, 
since $\overline{\tb}(\K)=-1$ and $w(\K)=0$.  In spite of this 
we have the following theorems:

\begin{thm}   \label{thm:torusknots}
Negative torus knots satisfy the {\rm UTP.}
\end{thm}

\begin{thm}  \label{thm:recursive}
If a knot type $\K$ satisfies the {\rm UTP,} then $(p,q)$\/{\rm -}\/cables 
$\K_{(p,q)}$ satisfies the {\rm UTP,} provided ${p\over q}< w(\K).$
\end{thm}

We sometimes refer to a slope ${p\over q}$ as ``sufficiently negative''
if ${p\over q}< w(\K).$\break\vskip-12pt\noindent Moreover, if ${p\over q}>w(\K)$ then we 
call the slope ``sufficiently positive''.

\begin{thm}  \label{thm:connectedsum}
If two knot types $\K_1$ and $\K_2$ satisfy the {\rm UTP,} then their connected 
sum $\K_1\# \K_2$ satisfies the {\rm UTP.}
\end{thm}

In Section~\ref{new} we give a more precise description and a proof of 
Theorem~\ref{negative cables} and in 
Section~\ref{section2} we prove Theorems~\ref{thm:torusknots} 
through \ref{thm:connectedsum} (the positive results on the UTP).

\Subsec{New phenomena}
While negative torus knots are well-behaved, positive torus knots 
are more unruly:  

\begin{thm}  \label{thm: positive torus}
There are positive torus knots that do not satisfy the {\rm UTP.}
\end{thm}

It is not too surprising that positive torus knots and negative torus 
knots have very different behavior --- recall that we also had to 
treat the positive and negative cases separately in the proof of the 
classification of Legendrian torus knots in \cite{EH1}.  A slight 
extension of Theorem~\ref{thm: positive torus} yields the following:

\begin{thm}  \label{thm: no destabilization}
There exist a knot type $\K$ and a Legendrian knot 
$L\in\mathcal{L}(\K)$ which does not admit any destabilization{\rm ,} yet 
satisfies $\tb(L) < \overline{\tb}(\K)$. 
\end{thm}

Although the phenomenon that appears in Theorem~\ref{thm: no 
destabilization} is rather common, we will specifically treat the case 
when $\K$ is a $(2,3)$-cable of a $(2,3)$-torus knot.  The same knot 
type $\K$ is also the example in the following theorem:

\begin{thm} \label{thm: transverse}
Let $\K$ be the $(2,3)$\/{\rm -}\/cable of the $(2,3)$\/{\rm -}\/torus knot. There is a
unique transverse knot in $\mathcal{T} (\K)$ for each self\/{\rm -}\/linking
number $n${\rm ,} where $n\leq 7$ is an odd integer $\not =3${\rm ,} and
exactly two transverse knots in $\mathcal{T} (\K)$ with self-linking
number $3$.  In particular{\rm ,} $\K$ is not transversely simple.
\end{thm} 

Here $\mathcal{T}(\K)$ is the set of transverse isotopy classes of $\K$.

Previously, Birman and Menasco \cite{BM} produced
non-transversely-simple knot types by exploiting an interesting
connection between transverse knots and closed braids.  It should be
noted that our theorem contradicts results of Menasco in
\cite{Menasco01}.  However, this discrepancy has led Menasco to find
subtle and interesting properties of cabled braids (see
\cite{Menasco?}).  The earlier work of Birman-Menasco \cite{BM} and our
Theorem~\ref{thm: transverse} both give negative answers to a
long-standing question of whether the self-linking number and the
topological type of a transverse knot determine the knot up to
contact isotopy.  The corresponding question for Legendrian knots, namely 
whether every topological knot type $\K$ is Legendrian simple, has 
been answered in the negative in the works of Chekanov \cite{Ch} and 
Eliashberg-Givental-Hofer \cite{EGH}.  Many other
non-Legendrian-simple knot types have been found since then (see for
example \cite{Ng}, \cite{EH2}). 

The theorem which bridges the Legendrian classification and the 
transverse classification is the following theorem from \cite{EH1}:

\begin{thm}   \label{stably}
Transverse simplicity is equivalent to {\em 
stable simplicity}{\rm ,} i.e.{\rm ,} any two $L_1,L_2\in \L(\K)$ with the same 
$\tb$ and $r$ become contact isotopic after some number of negative  
stabilizations. 
\end{thm}

The problem of finding a knot type which is not stably simple 
is much more difficult than the problem of finding a knot type which is not 
Legendrian simple, especially since the Chekanov-Eliashberg contact 
homology invariants vanish on stabilized knots.  Our technique for 
distinguishing stabilizations of Legendrian knots is to use the 
standard cut-and-paste contact topology techniques, and, in 
particular, the method of {\em state traversal}.

Theorems~\ref{thm: positive torus} and ~\ref{thm: no
  destabilization} will be proven in Section~\ref{section3} while
  Theorem~\ref{thm: transverse} will be proven in
  Section~\ref{sec:nts}. More specifically, the discussion in
  Section~\ref{sec:nts} provides a complete classification of
  Legendrian knots in the knot type of the
  $(2,3)$-cable of $(2,3)$-torus knot. 

\begin{figure}[ht]
$$
\epsfig{file=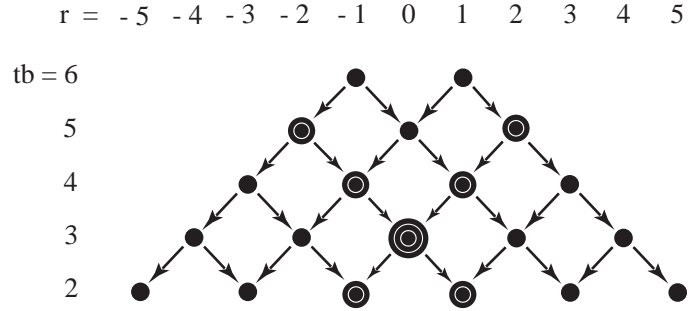}
$$
\caption{Classification of Legendrian $(2,3)$-cables of $(2,3)$-torus 
knots.  Concentric circles indicate multiplicities, i.e., the number of 
distinct isotopy classes with a given $r$ and $\tb$.}   
\label{iterate}  
\end{figure}

\begin{thm}
If $\K'$ is the $(2,3)$-cable of the $(2,3)$-torus knot{\rm ,} then 
$\L(\K')$ is classified as in Figure~{\rm \ref{iterate}.}  This entails 
the following\/{\rm :}\/
\begin{enumerate}
\item There exist exactly two maximal Thurston-Bennequin 
representatives $K_\pm\in \L(\K')$.  They satisfy $\tb(K_\pm)=6${\rm ,} 
$r(K_\pm)=\pm 1$. 
\item There exist exactly two non-destabilizable representatives 
$L_\pm\in \L(\K')$ which have non-maximal Thurston-Bennequin 
invariant. They satisfy $\tb(L_\pm)=5$ and $r(L_\pm)=\pm 2$. 
\item Every $L\in \L(\K')$ is a 
stabilization of one of $K_+${\rm ,} $K_-${\rm ,} $L_+${\rm ,} or $L_-$. 
\item $S_+(K_-)=S_-(K_+)${\rm ,}  $S_-(L_-)= S_-^2(K_-)${\rm ,} and 
$S_+(L_+)=S_+^2(K_+).$
\item $S_+^k(L_-)$ is not \/{\rm (}\/Legendrian\/{\rm )}\/ isotopic to $S_+^k S_-(K_-)$ 
and $S_-^k(L_+)$ is not isotopic to $S_-^k S_+(K_+)${\rm ,} for all positive 
integers $k$.  Also{\rm ,} $S_-^2(L_-)$ is not isotopic to $S_+^2(L_+)$.
\end{enumerate} 
\end{thm}

%%%%%%%%%%%%%%%%%%%%%%%%%%%%%%%%%%%%%%%%%%%%%%%%%%%%%%%
\section{Preliminaries}
%%%%%%%%%%%%%%%%%%%%%%%%%%%%%%%%%%%%%%%%%%%%%%%%%%%%%%%

Throughout this paper, a convex surface $\Sigma$ is either closed or
compact with Legendrian boundary, $\Gamma_\Sigma$ is the dividing set
of $\Sigma$, and $\#\Gamma_\Sigma$ is the number of connected
components of $\Gamma_\Sigma$.

\Subsec{Framings}
For convenience we relate the framing conventions that are used 
throughout the paper.  In what follows, $X\setminus Y$ will denote the 
metric closure of the complement of $Y$ in $X$.
  
Let $\K$ be a topological knot type and $\K_{(p,q)}$ be its 
$(p,q)$-cable.  Let $N(\K)$ be a solid torus which represents $\K$. 
Suppose $K_{(p,q)}\in\K_{(p,q)}$ sits on $\bdry N(\K)$.  Take an 
oriented annulus $A$ with boundary on $\bdry N(\K_{(p,q)})$ so that 
$(\bdry N(\K_{(p,q)}))\setminus A$ consists of two disjoint annuli 
$\Sigma_1$, $\Sigma_2$ and $A\cup \Sigma_i$, $i=1,2$, is isotopic to 
$\bdry N(\K)$.  We define the following coordinate systems, {\em i.e.}, 
identifications of tori with $\R^2/\Z^2$.

\begin{enumerate}
\item $\mathcal{C}_{\K}$, the coordinate system on $\bdry N(\K)$ where 
the (well-defined) longitude has slope $\infty$ and the meridian has 
slope $0$.
\item
$\mathcal{C}'_{\K}$, the coordinate system on $\bdry 
N(\K_{(p,q)})$ where the meridian has slope $0$ and slope $\infty$ 
is given by $A\cap \bdry N(\K_{(p,q)})$.
\end{enumerate}

We now explain how to relate the framings $\mathcal{C}'_{\K}$ and 
$\mathcal{C}_{\K_{(p,q)}}$ for $\bdry N(\K_{(p,q)})$.  Suppose 
$K_{(p,q)}\in \K_{(p,q)}$ is contained in $\bdry N(\K)$.  Then the 
Seifert surface $\Sigma(K_{(p,q)})$ is obtained by taking $p$ parallel 
copies of the meridional disk of $N(\K)$ (whose boundary we 
assume are $p$ parallel closed curves on $\bdry N(\K)$ of slope 0) and 
$q$ parallel copies of the Seifert surface for $\K$ (whose 
boundary we assume are $q$ parallel closed curves on $\bdry N(\K)$ of 
slope $\infty$), and attaching a band at each intersection between the 
slope 0 and slope $\infty$ closed curves for a total of $|pq|$ bands.  
Therefore, the framing coming from $\mathcal{C}'_{\K}$ and the framing 
coming from $\mathcal{C}_{\K_{(p,q)}}$ differ by $pq$; more precisely, 
if $L_{(p,q)}\in \mathcal{L}(\K_{(p,q)})$ and 
$t(L_{(p,q)},\mathcal{F})$ is the {\em twisting number} with respect 
to the framing $\mathcal{F}$ (or the {\em Thurston-Bennequin 
invariant} with respect to $\mathcal{F}$), then: 

\begin{equation} \label{eqn1}
t(L_{(p,q)},\mathcal{C}'_{\K}) +pq = 
t(L_{(p,q)},\mathcal{C}_{\K_{(p,q)}})= \tb(L_{(p,q)}). 
\end{equation} 

Let us also define the {\em maximal twisting number} of $\K$ with 
respect to $\mathcal{F}$ to be:
$$\overline{t}(\K,\mathcal{F}) = \max_{L\in \mathcal{L}(\K)} 
t(L,\mathcal{F}).$$

\Subsec{Computations of $\tb$ and $r$}
Suppose $L_{(p,q)}\in \mathcal{L}(\K_{(p,q)})$ is contained in $\bdry 
N(\K)$, which we assume to be convex.  We compute $\tb(L_{(p,q)})$ for 
two typical situations; the proof is an immediate consequence of 
equation~\ref{eqn1}.

\begin{lemma} \label{tbcomp}   $\mbox{ }$
\begin{enumerate}
\item Suppose $L_{(p,q)}$ is a Legendrian divide and 
${\rm slope} (\Gamma_{\bdry N(\K)})={q\over p}$.  Then 
$\tb(L_{(p,q)})=pq$. 
\item Suppose $L_{(p,q)}$ is a Legendrian ruling curve 
and ${\rm slope}  (\Gamma_{\bdry N(\K)})={q'\over p'}$. Then  
$\tb(L_{(p,q)})=pq- |pq'-qp'|.$ 
\end{enumerate}
\end{lemma}

Next we explain how to compute the rotation number $r(L_{(p,q)})$.  

\begin{lemma}  \label{rcomp}
Let $D$ be a convex meridional disk of $N(\K)$ with Legendrian 
boundary on a contact-isotopic copy of the convex surface $\bdry N(\K)${\rm ,}
 and 
let $\Sigma(L)$ be a convex Seifert surface with Legendrian boundary 
$L\in\mathcal{L}( \K)$ which is contained in a contact-isotopic copy of 
$\bdry N(\K)$.  \/{\rm (}\/Here the isotopic copies of $\bdry N(\K)$ are copies 
inside an $I$-invariant neighborhood of $\bdry N(\K)${\rm ,} obtained by 
applying the Flexibility Theorem to $\bdry N(\K)$.{\rm )}
Then $$ r(L_{(p,q)})= p \cdot r(\bdry D) + q \cdot r(\bdry \Sigma(K)).$$ 
\end{lemma}

\Proof 
Take $p$ parallel copies $D_1,\dots,D_p$ of 
$D$ and $q$ parallel copies $\Sigma(K)_1,\break\dots,\Sigma(K)_q$ of 
$\Sigma(K)$.  The key point is to use the Legendrian realization 
principle \cite{H1} simultaneously on $\bdry D_i$, $i=1,\dots,p$, and 
$\bdry \Sigma(K)_j$, $j=1,\dots,q$.
Provided ${\rm slope} (\Gamma_{\bdry N(\K)})\not 
=\infty$, the Legendrian realization principle allows us to 
perturb $\bdry N(\K)$ so that (i) $(\bigcup_{i=1,\dots,p} \bdry 
D_i)\cup (\bigcup_{j=1,\dots,q} \bdry \Sigma(K)_j)$ is a Legendrian 
graph in $\bdry N(\K)$ and (ii) each $\bdry D_i$ and 
$\bdry\Sigma(K)_j$ intersects $\Gamma_{\bdry N(\K)}$ {\em 
efficiently}, i.e., in a manner which minimizes the geometric 
intersection number.   (The version of Legendrian realization 
described in \cite{H1} is stated only for multicurves, but the proof 
for nonisolating graphs is identical.)  Now, suppose $L'_{(p,q)}\in 
\mathcal{L}(\K_{(p,q)})$ and its Seifert surface $\Sigma(L'_{(p,q)})$
are constructed by resolving the intersections of 
$(\bigcup_{i=1,\dots,p} \bdry D_i)\cup (\bigcup_{j=1,\dots,q} \bdry 
\Sigma(K)_j)$.  Recalling that the rotation number is a homological 
quantity (a relative half-Euler class) \cite{H1}, we readily compute 
that 
$$r(L'_{(p,q)}) = p \cdot r(\bdry D) + q \cdot r(\bdry \Sigma(K)).$$ 
(For more details on a similar computation, see \cite{EH1}.)  Finally, 
$L_{(p,q)}$ is obtained from $L'_{(p,q)}$ by resolving the inefficient 
intersections between $L'_{(p,q)}$ and $\Gamma_{\bdry N(\K)}$.  Since 
$\bdry N(\K)$ is a torus and $\Gamma_{\bdry N(\K)}$ consists of two 
parallel essential curves, the inefficient intersections come in 
pairs, and have no net effect on the rotation number computation.  
This proves the lemma. 
\hfill\qed

%%%%%%%%%%%%%%%%%%%%%%%%%%%%%%%%%%%%%%%%%%%%%%%%%%%%%%%
\section{From the UTP to classification}  \label{new}
%%%%%%%%%%%%%%%%%%%%%%%%%%%%%%%%%%%%%%%%%%%%%%%%%%%%%%%

In this section we use Theorem~\ref{thm:recursive} to give a complete
classification of $\L(\K_{(p,q)})$, provided $\L(\K)$ is classified,
$\K$ satisfies the UTP, and $\K$ is Legendrian simple.  In summary, we
show: 

\begin{thm}
If $\K$ is Legendrian simple and satisfies the {\rm UTP,} then all its
cables are Legendrian simple.
\end{thm}

The form of classification for Legendrian knots in the cabled knot
types depends on whether or not the cabling slope ${p\over q}$ is
greater or less than $w(\K)$.  The precise classification for
sufficiently positive slopes is given in Theorem~\ref{poscableclass},
while the classification for sufficiently negative slopes is given in
Theorem~\ref{negcableclass}. 
   
In particular, these results yield a complete classification of
Legendrian iterated torus knots, provided each iteration is
sufficiently negative (so that the UTP is preserved).  We follow the
strategy for classifying Legendrian knots as outlined in \cite{EH1}.

Suppose $\K$ satisfies the UTP and is Legendrian simple.  By the UTP, 
every Legendrian knot $L\in \mathcal{L}(\K)$ with 
$\tb(L)<\overline{\tb}(\K)$ can be destabilized to one realizing 
$\overline{\tb}(\K)$.  The Bennequin inequality \cite{Be} gives 
bounds on the rotation number; hence there are only finitely 
many distinct $L\in \mathcal{L}(\K)$, say $L_0,\dots, L_n$, which have 
$\tb(L_i)=\overline{\tb}(\K)$, $i=0,\dots,n$.  Write $r_i=r(L_i)$, and 
assume $r_0<r_1<\dots <r_n$.  By symmetry, $r_i=-r_{n-i}$. (This is
easiest to see in the front projection by rotating about the $x$-axis,
if the contact form is $dz-ydx$.)  Now, every time a Legendrian knot $L$ 
is stabilized by adding a zigzag, its $\tb$ decreases by $1$ and its 
$r$ either increases by $1$ (positive stabilization $S_+(L)$) or 
decreases by $1$ (negative stabilization $S_-(L)$).  Hence the image 
of $\mathcal{L}(K)$ under the map $(r,\tb)$ looks like a mountain 
range, where the peaks are all of the same height 
$\overline{\tb}(\K)$, situated at $r_0,\dots,r_n$.  The slope to the 
left of the peak is $+1$ and the slope to the right is $-1$, and the 
slope either continues indefinitely or hits a slope of the opposite 
sign descending from an adjacent peak to create a valley.  See 
Figure~2. \setcounter{figure}{1}
\begin{figure}[ht]
      $$
\epsfig{file=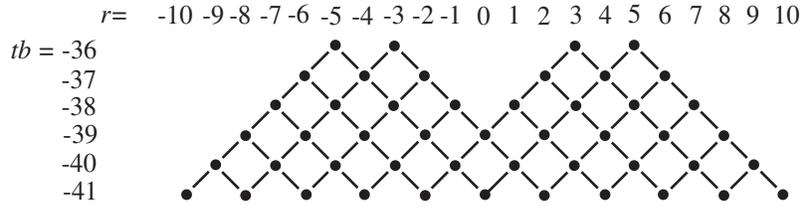}
$$
        \caption{The $(r,\tb)$-mountain range for the $(-9,4)$-torus knot.}
        \label{fig:geog}
\end{figure}

The following notation will be useful in the next few results.
Given two slopes $s=\frac{r}{t}$ and $s'=\frac{r'}{t'}$ on a torus $T$
with $r,t$ relatively prime and $r',t'$ relatively prime, we denote:
\[s\bullet s'=rt'-tr'.\]
This quantity is the minimal number of intersections between two
curves of slope $s$ and $s'$ on $T.$

\begin{thm}\label{poscableclass}
Suppose $\K$ is Legendrian simple and satisfies the {\rm UTP}.  If $p,q$ are 
relatively prime integers with $\frac{p}{q}> w(\K)${\rm ,} then
$\K_{(p,q)}$ is also Legendrian simple.  Moreover{\rm ,}
$$\overline{tb}(\K_{(p,q)})=pq-\left|w(\K)\bullet
\frac{p}{q}\right|,$$ 
and the set of rotation numbers realized by $\{L\in \L(\K_{(p,q)})|
tb(L)=\overline{tb}(\K_{(p,q)})\}$ is 
\[\{q\cdot r(L)| L\in \L(\K), tb(L)=w(\K)\}.\]
\end{thm}

This theorem is established through the following three lemmas.

\begin{lemma}
Under the hypotheses of Theorem~{\rm \ref{poscableclass},}
$\overline{tb}(\K_{(p,q)})=pq-|w(\K)\bullet\frac{p}{q}|$ and any
Legendrian knot $L\in\L(\K_{(p,q)})$ with
$tb(L)<\overline{tb}(\K_{(p,q)})$ destabilizes. 
\end{lemma}

\Proof 
We first claim that $t(L, \CC'_K)<0$ for any $L\in \L(\K_{(p,q)})$.  If
not, there exists a Legendrian knot $L'\in \L(\K_{(p,q)})$ with $t(L',
\CC'_K)=0$.  Let $S$ be a solid torus representing $\K$ such that
$L'\subset \partial S$ (as a Legendrian divide) and the boundary torus
$\partial S$ is convex.  Then ${\rm slope} (\Gamma_{\bdry
  S})=\frac{q}{p}$ when measured with respect to $\CC_K.$  However,
since $\frac{p}{q}>w(\K)$, this contradicts the UTP.   

Since $t(L, \CC'_\K)<0$, there exists an $S$ so that $L\subset \bdry
S$ and $\partial S$ is convex.  Let $s$ be the slope of
$\Gamma_{\partial S}$.  Then we have the following inequality:
$$\left|\frac{1}{s}\bullet \frac{p}{q}\right|\geq \left|w(\K)\bullet
\frac{p}{q}\right|,$$ 
with equality if and only if $\frac{1}{s}=w(\K).$  To see this, use an
oriented diffeomorphism of the torus $\partial S$ that sends slope $0$
to $0$ and slope $\frac{1}{w(\K)}=\frac{1}{\overline{tb}(\K)}$ to
$\infty$ (this forces $-\infty \leq s < 0$ and
$\frac{q}{p}>0$), and compute determinants. (Alternatively, this follows
from observing that there is an edge from $0$ to $\frac{1}{w(\K)}$ in the
Farey tessellation, and $\frac{1}{s}\in (-\infty,w(\K)]$, whereas ${p\over
q}\in(w(\K),\infty)$.)  Thus $t(L,\CC'_\K)\leq-|w(\K)\bullet\frac{p}{q}|$ 
for all $L\in\L(\K_{(p,q)})$.  But now, if $S$ is a solid torus
representing $\K$ of maximal thickness, then a Legendrian ruling curve
on $\bdry S$ easily realizes the equality.  Converting from $\CC'_\K$
to $\CC_\K$, we obtain $\overline{tb}(\K_{(p,q)})=
pq-|w(\K)\bullet\frac{p}{q}|.$   

Now consider a Legendrian knot $L\in \L(\K_{(p,q)})$ with
$tb(L)<\overline{tb}(\K_{(p,q)})$.  Placing $L$ on a convex surface
$\bdry S$, if the intersection between $L$ and $\Gamma_{\bdry S}$ is
not efficient (i.e., does not realize the geometric intersection
number), then there exists a bypass which allows us to destabilize
$L$.  Otherwise $L$ is a Legendrian ruling curve on $\bdry S$ with
$\frac{1}{s}\not=w(\K)$.  Now, since $\K$ satisfies the UTP, there is
a solid torus $S'$ with $S\subset S'$, where $\partial S'$ is convex
and ${\rm slope} (\Gamma_{\bdry S'})= \frac{1}{w(\K)}$.  By comparing
with a Legendrian ruling curve of slope ${q\over p}$, i.e., taking a
convex annulus $A=L\times [0,1]$ in $\partial S\times [0,1]=
S'\setminus S$ and using the Imbalance Principle, we may easily find a
bypass for $L$.  Therefore, if $t(L, \CC'_\K)< 
-|w(\K)\bullet \frac{p}{q}|$, then we may destabilize $L.$ 
\Endproof\vskip4pt  

\begin{lemma}\label{protation}
Under the hypotheses of Theorem~{\rm \ref{poscableclass},} Legendrian knots
with maximal $tb$ in $\L(\K_{(p,q)})$ are determined by their rotation
number.  Moreover{\rm ,} the rotation numbers associated to maximal $tb$
Legendrian knots in $\L(\K_{(p,q)})$ are \[\{q\cdot r(L)| L\in \L(\K),
tb(L)=w(\K)\}.\] 
\end{lemma}
\vskip9pt

\Proof 
Given a Legendrian knot $L\in\L(\K_{(p,q)})$ with maximal $tb$, there 
exists a solid torus $S$ with convex boundary, where 
${\rm slope} (\Gamma_{\partial S})= \frac{1}{w(\K)}$ and $L$\break\vskip-12pt\noindent  is a
Legendrian ruling curve on $\partial S.$  The torus $S$ is a standard
neighborhood of a Legendrian knot $K$ in $\L(\K).$ From
Lemma~\ref{rcomp} one sees that \[r(L)=q\cdot r(K).\] Thus the
rotation number of $L$ determines the rotation number of $K.$   

If $L$ and $L'$ are two Legendrian knots in $\L(\K_{(p,q)})$ with
maximal $tb$, then we have the associated solid tori $S$ and $S'$ and
Legendrian knots $K$ and $K'$ as above.  If $L$ and $L'$ have the same
rotation numbers then so do $K$ and $K'.$ Since $\K$ is Legendrian
simple, $K$ and $K'$ are Legendrian isotopic. Thus we may assume that
$K$ and $K'$ are the same Legendrian knot and that $S$ and $S'$ are
two standard neighborhoods of $K=K'.$ Inside $S\cap S'$ we can find
another standard neighborhood $S''$ of $K=K'$ with convex boundary
having dividing slope $\frac{1}{w(\K)}$ and ruling slope
$\frac{q}{p}.$   The sets 
$S\setminus S''$ and $S'\setminus S''$ are both diffeomorphic to
$T^2\times [0,1]$ and have $[0,1]$-invariant contact structures. Thus
we can assume that $L$ and $L'$ are both ruling curves on $\partial
S''.$ One may now use the other ruling curves on $\partial S''$ to
Legendrian isotop $L$ to $L'.$   
\hfill\qed

\begin{lemma}
Under the hypotheses of Theorem~{\rm \ref{poscableclass},} Legendrian knots
in $\L(\K_{(p,q)})$ are determined by their Thurston-Bennequin
invariant and rotation number. 
\end{lemma}

\Proof 
Here one simply needs to see that there is a unique Legendrian knot in
 the valleys of the $(r,\tb)$-mountain range; that is, if $L$ and $L'$
 are maximal $tb$ Legendrian knots in $\L(\K_{(p,q)})$ and
 $r(L)=r(L')+2qn$ (note the difference in their rotation numbers must
 be even and a multiple of $q$) then $S_+^{qn}(L')=S_-^{qn}(L).$ To
 this end, let $K$ and $K'$ be the Legendrian knots in $\L(\K)$
 associated to $L$ and $L'$ as in the proof of the previous lemma. The
 knots $K$ and $K'$ have maximal $tb$ and $r(K)=r(K')+2n.$ Since $\K$
 is Legendrian simple we know $S_+^n(K')=S_-^n(K).$ Using the fact
 that $S_-^q(L)$ sits on a standard neighborhood of $S_-(K)$ (and the
 corresponding fact for $K'$ and $L'$) it easily follows that
 $S_+^{qn}(L')=S_-^{qn}(L).$ 
\Endproof\vskip4pt  

We now focus our attention on sufficiently negative cablings of a
knot type~$\K$.

\begin{thm}\label{negcableclass}
Suppose $\K$ is Legendrian simple and satisfies the {\rm UTP}.  If $p,q$ are
relatively prime integers with $q>0$ and $\frac{p}{q}< w(\K)${\rm ,} then
$\K_{(p,q)}$ is also Legendrian simple.  Moreover
$\overline{tb}(\K_{(p,q)})=pq$ and the set of rotation numbers realized by 
$\{L\in \L(\K_{(p,q)})| tb(L)=\overline{tb}(\K_{(p,q)})\}$ is 
\[\{\pm (p + q(n+r(L))) \mbox{ } | \mbox{ } L\in \L(\K), \tb(L)=-n \},\]
where $n$ is the integer that satisfies
$$-n-1 < {p\over q}< -n.$$
\end{thm}

We begin with two lemmas.

\begin{lemma}\label{ndes}
Under the hypotheses of Theorem~{\rm \ref{negcableclass},} every
$L_{(p,q)}\in \L(\K_{(p,q)})$ with  
$\tb(L_{(p,q)})< \overline{\tb}(\K_{(p,q)})$ can be destabilized and
$$\overline{\tb}(\K_{(p,q)})=pq.$$
\end{lemma}

\Proof 
By Theorem~\ref{thm:recursive}, $\K_{(p,q)}$ also satisfies the UTP.
Therefore\break every $L_{(p,q)}\in \L(\K_{(p,q)})$ with 
$\tb(L_{(p,q)}) < \overline{\tb}(\K_{(p,q)})$ can be destabilized to a\break
Legendrian knot realizing $\overline{\tb}(\K_{(p,q)})$.  Moreover, 
since ${p\over q}$ is sufficiently negative, there exist $L_{(p,q)}\in 
\L(\K_{(p,q)})$ with $\tb(L_{(p,q)})=pq$, which appear as Legendrian
divides on a convex torus $\bdry N(\K)$. By Lemma~\ref{tbcomp} we have 
$\overline{\tb}(\K_{(p,q)})\geq pq$.  Equality (the hard part) follows from
Claim~\ref{cabletbbound} below. 
\hfill\qed

\begin{lemma}\label{rotation}
Under the hypotheses of Theorem~{\rm \ref{poscableclass},} Legendrian knots
with maximal $tb$ in $\L(\K_{(p,q)})$ are determined by their rotation
number. Moreover{\rm ,} the set of rotation numbers attained by $\{L_{(p,q)} \in 
\L(\K_{(p,q)}) \mbox{ } | \mbox{ } \tb(L_{(p,q)})=pq\}$ is $$\{\pm (p 
+ q(n+r(L))) \mbox{ } | \mbox{ } L\in \L(\K), \tb(L)=-n \}.$$ 
\end{lemma}

Another way of stating the range of rotation numbers (and seeing where
they come from) in Lemma~\ref{rotation} is as follows:  To each 
$L\in \L(\K)$, there correspond two elements $L^\pm\in \L(\K_{(p,q)})$ 
with $\tb(L^\pm)=pq$ and $r(L^\pm)=q\cdot r(L) \pm s$, where $s$ is 
the remainder $s=-p-qn>0$.  $L^\pm$ is obtained by removing a standard 
neighborhood of $N(S_\pm(L))$ from $N(L)$, and considering a Legendrian 
divide on a torus with ${\rm slope} (\Gamma)={q\over p}$ inside 
$T^2\times[1,2]= N(L)\setminus N(S_\pm(L))$.   

\Proof 
The proof that Legendrian knots with maximal $tb$ in $\L(\K_{(p,q)})$
are determined by their rotation numbers
is similar to the proof of Lemma~\ref{protation} (also see \cite{EH1}).

The range of rotation numbers follows from Lemma~\ref{rcomp} as well
as some considerations of tight contact structures on thickened tori.
First let $T_{1.5}=\bdry N(\K)$ which contains $L_{(p,q)}$ with 
$\tb(L_{(p,q)})=pq$.  We will use the coordinate system 
$\mathcal{C}_\K$.  Then there exists a thickened torus 
$T^2\times[1,2]$ with convex boundary, where $T^2\times[1,1.5]\subset 
N(\K)$, ${\rm slope} (\Gamma_{T_1})=-{1\over n+1}$,
${\rm slope} (\Gamma_{T_{1.5}})={q\over p}$, and
${\rm slope} (\Gamma_{T_2})=-{1\over n}$.  Here we write 
$T_i=T^2\times\{i\}$.  Observe that $T^2\times[1,2]$ is a {\em basic 
slice} in the sense of \cite{H1}, since the shortest integral vectors 
$(-n,1)$ and $(-(n+1),1)$ form an integral basis for $\Z^2$.  This 
means that the tight contact structure must be one of two 
possibilities, distinguished by the {\em relative half-Euler class}
$e(\xi)$.  (It is called the ``relative Euler class'' in \cite{H1},
but ``relative half-Euler class'' is more appropriate.)  Their
Poincar\'e duals are given by  $PD(e(\xi))=\pm ((-n,1) -(-(n+1),1))=
\pm (1,0)$.  Now, by  
the universal tightness of $T^2\times[1,2]$, it follows from the 
classification of \cite{Gi2}, \cite{H1} that:
\begin{enumerate}
\item  either $PD(e(\xi),T^2\times[1,1.5])=(p,q)- (-n-1,1)$ and 
$PD(e(\xi),T^2\times[1.5,2])=(-n,1) -(p,q)$,
\item or $PD(e(\xi),T^2\times[1,1.5])=-(p,q)+ (-n-1,1)$ and 
$PD(e(\xi),T^2\times[1.5,2])=-(-n,1) + (p,q).$
\end{enumerate}

In view of Lemma~\ref{rcomp}, we want to compute (i) $r(\bdry D)$, 
where $D$ is a convex meridional disk for $N(\K)$ with Legendrian 
boundary on $T_{1.5}=\bdry N(\K)$, and (ii) $r(\bdry \Sigma)$, where 
$\Sigma$ is a convex Seifert surface for a Legendrian ruling curve 
$\infty$ on $T_{1.5}$.  Write $D=D'\cup A$, where $D'$ is a meridional 
disk with {\em efficient} Legendrian boundary for $N(\K)\setminus 
(T^2\times[1,1.5])$, and $A\subset T^2\times[1,1.5]$.   (An {\em 
efficient} closed curve on a convex surface intersects the dividing 
set $\Gamma$ minimally.)   Also write $\Sigma=\Sigma'\cup B$, where 
$B\subset T^2\times[1.5,2]$ and $\Sigma'\subset S^3\setminus 
(T^2\times[1,2])$ has efficient Legendrian boundary $L$ on $T_2$.

By additivity, 
$$r(\bdry \Sigma)= r(L) + \chi (B_+) -\chi (B_-) = 
r(L) + \langle e(\xi), B\rangle. $$  
Here  $S_+$ (resp.\ $S_-$) denotes the positive (resp.\ negative)
region of a convex surface $S$, divided by $\Gamma_S$.  Similarly, 
$$r(\bdry D) = r(\bdry D') + \langle e(\xi), A\rangle = \langle 
e(\xi),A\rangle.$$
Therefore,  either $r(\bdry \Sigma) = r(L) + p+n$ and $r(\bdry 
D)=-q+1$, or $r(\bdry \Sigma) = r(L) - p-n$ and $r(\bdry D)=q-1$.  In
the former case, 
$$r(L_{(p,q)})=p(-q+1) + q(r(L) +p+n)=p +q(r(L)+n).$$
In the latter case, we have $r(L_{(p,q)})= -p +q(r (L)-n)$ and we use 
the fact that $\{r(L)\mbox{ }|\mbox{ } L\in \L(\K), \tb(L)=-n\}$ is 
invariant under the map $r\mapsto -r$. 
\Endproof\vskip4pt  

{\it Proof of Theorem}~\ref{negcableclass}.
By Lemma~\ref{ndes}, every $L_{(p,q)}'\in 
\L(\K_{(p,q)})$ can be written as $S_+^{k_1} S_-^{k_2} (L_{(p,q)})$ 
for some $L_{(p,q)}$ with maximal $\tb$.  To complete the 
classification, we need to show that every $L_{(p,q)}'$ which is a 
``valley'' of the image of $(r,\tb)$ (i.e., $L_{(p,q)}'$ for which 
$(r(L_{(p,q)}')\pm 1,\tb(L_{(p,q)}')+1)$ is in the image of $(r,\tb)$ 
but $(r(L_{(p,q)}'), \tb(L_{(p,q)}')+2)$ is not) destabilizes to two 
maximal $\tb$ representatives $L_{(p,q)}^+$ and $L_{(p,q)}^-$ (the 
``peaks'').  Observe that there are two types of valleys:  
type (i) has a depth of $s=-p-qn$ and type (ii) has a depth 
of $kq -s$, $k\in \Z^+$.  

We start with valleys of type (i).  Such valleys occur when 
$r(L^-)= q\cdot r(L) -s$, $r(L^+)= q\cdot r(L)+s$, 
and $\tb(L^-)=\tb(L^+)=pq$.  It is clear that the valley between $L^-$ 
and $L^+$ corresponds to a Legendrian ruling curve of slope ${q\over 
p}$ on the boundary of the standard neighborhood $N(L)$ of $L$ with 
$\tb(L)=-n$.   By stabilizing $L$ in two ways, we see that any element 
$L_{(p,q)}'$ with $r(L_{(p,q)}')= q\cdot r(L)$ and $tb(L_{(p,q)}')= pq 
-s$ satisfies $L_{(p,q)}'= S_-^s(L^+)=S_+^s(L^-).$ 

Next we explain the valleys of type (ii) which have depth
$kq-s$, $k\in \Z^+$.  The peaks $L^-$ and $L^+$ 
correspond to ``adjacent'' $L,L'\in \L(\K)$ which have $\tb(L)=\tb(L') 
=-n$ and $r(L) < r(L')$, and such that there is no Legendrian $L''\in 
\L(\K)$ with $\tb(L'')=-n$ and $r(L)<r(L'') <r(L')$.  Hence $r(L^-)= 
q\cdot r(L) +s$ and $r(L^+)= q\cdot r(L') -s$.  The $k$ in the 
expression $kq-s$ above satisfies $r(L')-r(L)= 2k$.  The valley 
$L_{(p,q)}'$ with $\tb(L_{(p,q)}')= pq-(kq-s)$ and $r(L_{(p,q)}')= 
q\cdot r(L) +kq= q\cdot r(L')-kq$ occurs as a Legendrian ruling 
curve of slope ${q\over p}$ on the standard tubular neighborhood of 
$S_+^{k}(L) = S_-^{k}(L')$.  Therefore, $L_{(p,q)}'= 
S_+^{kq-s}(L^-)=S_-^{kq-s}(L^+)$.  This proves the Legendrian 
simplicity of $\L(\K_{(p,q)})$.
\hfill\qed

%%%%%%%%%%%%%%%%%%%%%%%%%%%%%%%%%%%%%%%%%%%%%%%%%%%%%%%
\section{Verification of uniform thickness}  \label{section2}
%%%%%%%%%%%%%%%%%%%%%%%%%%%%%%%%%%%%%%%%%%%%%%%%%%%%%%%

In this section we prove that many knot types satisfy the UTP.  
Let us begin with negative torus knots.

\demo{\scshape Theorem~\ref{thm:torusknots}}
{\it Negative torus knots satisfy the} UTP.

\Proof 
Let $\K$ be the unknot and $\K_{(p,q)}$ be its $(p,q)$-cable, {\em i.e.}, 
the  $(p,q)$-torus knot, with $pq<0$.  It was shown 
in \cite{EH1} that $\overline{\tb}(\K_{(p,q)})= pq$.  Unless indicated 
otherwise, we measure the slopes of tori isotopic to $\bdry 
N(\K_{(p,q)})$ with respect to $\mathcal{C}'_{\K}$.  Then 
$\overline{\tb}(\K_{(p,q)})=pq$ is equivalent to 
$\overline{t}(\K_{(p,q)})=\overline{t}(\K_{(p,q)},\mathcal{C}'_{\K})\break =0 
$.  In other words, the standard neighborhood of $L\in 
\mathcal{L}(\K_{(p,q)})$ satisfying $\tb(L)=pq$ has boundary slope 
$\infty$ with respect to $\mathcal{C}'_{\K}$.

We will first verify Condition~1 of the UTP, arguing by contradiction.  
(In fact, the argument that follows can be used to prove that 
$\overline{t}(\K_{(p,q)})=0$.)  Suppose there exists a solid torus 
$N=N(\K_{(p,q)})$ which has convex boundary with 
$s={\rm slope} (\Gamma_{\bdry N})>0$ and $\#\Gamma_{\bdry N}=2$.
After shrinking $N$ if necessary, we may assume that $s$ is a large positive
integer.  Next, using the Giroux  
Flexibility Theorem, $\bdry N$ can be isotoped into {\em standard 
form}, with Legendrian rulings of slope $\infty$.  Now let $A$ be a 
convex annulus with Legendrian boundary on $\bdry N$ and $A\times 
[-\varepsilon,\varepsilon]$ its invariant neighborhood.  Here $A$ is
chosen so that $R=N\cup (A\times[-\varepsilon,\varepsilon])$ is a
thickened torus whose boundary $\bdry R=T_1\cup T_2$ is parallel to
$\bdry N(\K)$.
Here, the relative positions of $T_1$ and $T_2$ are that if $T_2=\bdry
N(\K)$, then $T_1\subset N(\K)$.    

Let us now analyze the possible dividing sets for $A$. First,
$\bdry$-parallel dividing curves are easily eliminated.  Indeed, if
there is a $\bdry$-parallel arc, then we may attach the corresponding
bypass onto $\bdry N$ and increase $s$ to $\infty$, after isotopy.
This would imply excessive twisting inside $N$, and the contact
structure would be overtwisted.  Hence we may assume that $A$ is in 
standard form, with two parallel nonseparating arcs.  Now choose a
suitable identification $\bdry N(\K)\simeq \R^2/\Z^2$ so that the  
ruling curves of $A$ have slope $\infty$, ${\rm slope} (\Gamma_{T_1})=-s$ and
${\rm slope} (\Gamma_{T_2})=1$. (This is possible since a holonomy 
computation shows that $\Gamma_{T_1}$ is obtained from $\Gamma_{T_2}$
by performing $s+1$ right-handed Dehn twists.)

We briefly explain the classification of tight contact structures on
$R$ with the boundary condition ${\rm slope} (\Gamma_{T_1})=-s$,
${\rm slope} (\Gamma_{T_2})=1$,  
\hbox{$\#\Gamma_{T_1}=\#\Gamma_{T_2}=2$.} For more details, see \cite{H1}.
Corresponding to the slopes $-s$, $1$, are the shortest integer
vectors $(-1,s)$ and $(1,1)$.  Any tight contact structure on $R$ can
naturally be layered into {\em basic slices} $(T^2\times[1,1.5]) \cup
(T^2\times[1.5,2])$, where ${\rm slope} (\Gamma_{T_{1.5}})=\infty$
(corresponding to the shortest integer vector $(0,1)$) and
$\#\Gamma_{T_{1.5}}=2$.  There are two possibilities for each basic
slice --- the Poincar\'e duals of the relative half-Euler classes are
given by $\pm$ the difference of the shortest integer vectors
corresponding to the dividing sets on the boundary.  For
$T^2\times[1,1.5]$, the possible $PD(e(\xi))$ are $\pm$ of
$(0,1)-(-1,s)= (1,1-s)$; for $T^2\times[1.5,2]$, the possibilities are
$\pm$ of $(1,1)-(0,1)=(1,0)$.  Since $s\gg 1$, the four possible tight
contact structures on $R$ are given by $\pm (1,0) \pm (1,1-s)$.  Of
the four possibilities, two  of them are universally tight and two  of them are
virtually overtwisted.  The contact structure $\xi$ is universally
tight when there is {\em no mixing of sign}, {\em i.e.}, $PD(e(\xi))=
+ (1,0) + (1,1-s)$ or 
$-(1,0)-(1,1-s)$; when there is mixing of sign $+(1,0)-(1,1-s)$ or
$-(1,0)+(1,1-s)$, the contact structure is virtually overtwisted. 

To determine the half-Euler class, consider
$\Sigma=\gamma\times[-\varepsilon,\varepsilon]\subset
A\times[-\varepsilon,\varepsilon]$, where $\gamma$ is a Legendrian
ruling curve of slope $\infty$.  Since $\Sigma$ is
$[-\varepsilon,\varepsilon]$-invariant, $\langle 
e(\xi), \Sigma\rangle=\chi(\Sigma_+)-\chi(\Sigma_-)=0$, where $\chi$
is the Euler characteristic and $\Sigma_+$ (resp.\ $\Sigma_-$) is the
positive (resp.\ negative) part of $\Sigma\setminus \Gamma_\Sigma$.
Therefore, $PD(e(\xi))$ must be $\pm(0,s-1)$, implying a {\em mixture
  of sign}. 

Let us now recast the slopes of $\Gamma_{T_i}$ in terms of coordinates
$\mathcal{C}_{\K}$, where $\K$ is the unknot.  With respect to
$\mathcal{C}_{\K}$, ${\rm slope} (\Gamma_{T_{1.5}})= {q\over p}$,
where ${q\over p}$ is neither a negative integer nor the reciprocal
of one.  One of the consequences of the classification of tight
contact structures on solid tori in \cite{Gi2}, \cite{H1} is the following:
if $S$ is a convex torus in the standard tight contact
$(S^3,\xi_{std})$ which bounds solid tori on both sides, then the only
slopes for $\Gamma_S$ at which there can be a sign change are negative
integers or reciprocals of negative integers.  Therefore, we have a
contradiction, proving Condition~1. 

Next we prove Condition~2, keeping the same notation
as in the proof of Condition~1.  Suppose that $N=N(\K_{(p,q)})$ now
has boundary slope $s$, where $-\infty < s < 0$ and slopes are
measured with respect to $\mathcal{C}'_{\K}$. If $\Gamma_A$ has a
$\bdry$-parallel arc, then $s$ approaches $-\infty$ (in
terms of the Farey tessellation) when we attach a corresponding bypass
onto $N$.  Therefore, as usual, we may take $A$ to be in standard form
and $\Gamma_A$ to consist of parallel nonseparating dividing arcs.
Now observe that 
$\frac{q}{p}$ cannot lie between ${\rm slope} (\Gamma_{T_1})$ and
${\rm slope} (\Gamma_{T_2})$, where the slopes are measured with
respect to $\mathcal{C}_\K$. This implies that there are no convex
tori in $R$ which are isotopic to $T_i$ and have slope
$\frac{q}{p}$. In the complement $S^3\setminus R$, there is a convex
torus isotopic to $T_i$ with slope $\frac{q}{p}$.  
Using this, we readily find a thickening of $N$ to have 
slope $\infty$, measured with respect to $\mathcal{C}'_{\K}$.  

Once we thicken $N$ to have boundary slope $\infty$, there is one last 
thing to ensure, namely that $\#\Gamma_{\bdry N}=2$; in other words, 
we want $N$ to be the standard neighborhood of a Legendrian curve with 
twisting number $0$ with respect to~$\mathcal{C}'_\K$.

\begin{claim}   \label{reduction of number}
Any solid torus $N$ with convex boundary{\rm ,} ${\rm slope} (\Gamma_{\bdry 
N})=\infty${\rm ,} and $\#\Gamma_{\bdry N}=2n$, $n>1${\rm ,} extends to a solid torus 
$\overline{N}$ with convex boundary{\rm ,} slope $\infty$, and $\#\Gamma_{\bdry 
\overline{N}}=2$. 
\end{claim}

\Proof 
There exists a thickened torus $R$ with $\bdry R=T_2-T_1$, where
$N\subset R$, the $T_i$, $i=1,2$, bound solid tori on both sides, and
${\rm slope} (\Gamma_{T_i})={q\over p}$ with respect to
$\mathcal{C}_\K$.  By shrinking $N$ somewhat, we may take $R\setminus
N$ to be a
pair-of-pants $\Sigma_0$ times $S^1$.  Since there is
twisting on both sides of the exterior of $R$, we may also arrange that
$\#\Gamma_{T_i}=2$.  Moreover, as $\Gamma_{\bdry (R\setminus N)}$ is
parallel to the $S^1$-fibers, the tight contact structure on
$R\setminus N$ is necessarily {\em  
vertical}, i.e., isotopic to an $S^1$-invariant contact structure, 
after appropriately modifying the boundary to be Legendrian-ruled.  
(See \cite{H2} for a proof.) 

The data for this tight contact structure are encoded in
$\Gamma_{\Sigma_0}$.  (Here we are assuming without loss of generality
that $\Sigma_0$ is convex with Legendrian boundary.)  Let $\bdry
\Sigma_0= \gamma\sqcup \gamma_1\sqcup \gamma_2$, where
$\gamma_i=\Sigma_0\cap T_i$ and $\gamma=\bdry N\cap \Sigma_0$.  There
are $2n$ endpoints of $\Gamma_{\Sigma_0}$ on $\gamma$, and $2$ for
each of $\gamma_i$.   If there is an arc between $\gamma_1$ and
$\gamma_2$, then an imbalance occurs and there is necessarily a
$\bdry$-parallel arc along $\gamma$.  This would allow a thickening of
$N$ to one whose boundary has fewer dividing curves.  

The situation from which we have no immediate escape is when all the
arcs from $\gamma_i$ go to $\gamma$, and the extra endpoints along
$\gamma$ connect up without creating $\bdry$-parallel arcs.  We need
to look externally (i.e., outside of $R$) to obtain the desired bypass.
The key features we take advantage of are: 
\begin{enumerate}
\item There is twisting on both sides of the exterior of $R$.
\item There is no mixing of sign about $R$.
\end{enumerate} 
One of the (nontrivial) bypasses found along $T_1$ and $T_2$ therefore
can be extended into $R$ to give a bypass to reduce $\#\Gamma_{\bdry N}$. 
\Endproof\vskip4pt  

This completes the proof of Theorem~\ref{thm:torusknots}.
\Endproof\vskip4pt  

Recall a fraction ${p\over q}$ is {\em sufficiently negative} if 
$${p\over q} < w(\K).$$
(Observe that ${p\over q}$ is the reciprocal of the slope of a curve
$\bdry N$ corresponding to $(p,q)$.)   

\demo{\scshape Theorem~\ref{thm:recursive}}
{\it If a knot type $\K$ satisfies the {\rm UTP,} then $(p,q)$-cables
$\K_{(p,q)}$ satisfies the {\rm UTP,} provided ${p\over q}$ is sufficiently
negative. }
\Enddemo

Let $\K$ be a knot type that satisfies the UTP.  We write $N=N(\K)$
and $N_{(p,q)}=N(\K_{(p,q)})$.  The coordinates for $\bdry N$ and $\bdry 
N_{(p,q)}$ will be $\mathcal{C}_\K$ and $\mathcal{C}'_\K$, 
respectively.  The proof of Theorem~\ref{thm:recursive} is virtually
identical to that of Theorem~\ref{thm:torusknots}.

\Proof 
We prove that the contact width $w(\K_{(p,q)},\mathcal{C}'_\K)$, 
measured with respect to $\mathcal{C}'_\K$, and 
$\overline{t}(\K_{(p,q)},\mathcal{C}'_\K)$ both equal $0$, and that 
any $N_{(p,q)}$ with convex boundary can be thickened to a standard 
neighborhood of a Legendrian knot with 
$t(L_{(p,q)},\mathcal{C}'_\K)=0$. 

It is easy to see that $t(L_{(p,q)}, \mathcal{C}'_\K)=0$ can be attained: 
Since ${p\over q}$ is sufficiently negative, inside any $N$ (with 
convex boundary) of maximal thickness there exists a Legendrian
representative $L_{(p,q)}\in \L(\K_{(p,q)})$ of twisting number
$t(L_{(p,q)})=0$, which appears as a Legendrian divide on a convex
torus parallel to $\bdry N$.

Suppose $N_{(p,q)}$ has convex boundary and
${\rm slope} (\Gamma_{\bdry N_{(p,q)}})=s$.  As before, arrange the
characteristic foliation on $\bdry N_{(p,q)}$ to be in standard form
with Legendrian rulings of slope $\infty$, and consider the convex
annulus $A$ with Legendrian boundary on $\bdry N_{(p,q)}$, where the
thickening $R$ of $N_{(p,q)}\cup A$ is a thickened torus whose
boundary $\bdry R=T_1\cup T_2$ is isotopic to $\bdry N$.   We assume that 
$\Gamma_A$ consists of parallel nonseparating arcs, since otherwise we 
can further thicken $N_{(p,q)}$ by attaching the bypass corresponding 
to a $\bdry$-parallel arc. 

Now let $N$ be a maximally thickened solid torus which contains $R$, 
where the thickness is measured in terms of the {\em contact width}.  

\begin{claim}\label{cabletbbound}
$w(\K_{(p,q)},\mathcal{C}'_\K)=\overline{t}(\K_{(p,q)},\mathcal{C}'_\K
)=0.$
\end{claim}

\Proof 
If $s>0$, then by shrinking the solid torus 
$N_{(p,q)}$, we may take $s$ to be a large positive integer 
and $\#\Gamma_{\bdry N_{(p,q)}}=2$.  Then, as in the proof of 
Theorem~\ref{thm:torusknots}, (i) inside $R$ there exists a convex 
torus parallel to $T_i$ with slope ${q\over p}$ (with respect to 
$\mathcal{C}_\K$), (ii) the tight contact structure on $R$ must have 
mixing of sign, and (iii) this mixing of sign cannot happen inside the 
maximally thickened torus $N$.    This contradicts 
${\rm slope} (\Gamma_{\bdry N_{(p,q)}})=s>0$. 
\hfill\qed

\begin{claim}
Every $N_{(p,q)}$ can be thickened to a standard neighborhood of a 
Legendrian knot $L_{(p,q)}$ with $t(L_{(p,q)})=0$.
\end{claim}

\Proof 
If $-\infty < s <0$, then there cannot be any convex tori in 
$R$ isotopic to $T_i$ and with slope $\infty$.  Hence there is a 
convex torus parallel to $T_i$ with slope $\infty$ and $\#\Gamma=2$ 
{\em outside of} $R$.  By an application of the Imbalance Principle, 
we can thicken $N$ to have slope $\infty$.   The proof of the 
reduction to $\#\Gamma_{\bdry N}=2$ is identical to the proof of 
Claim~\ref{reduction of number} --- the key point is that there is 
twisting on both sides of $N\setminus R$. 
\Endproof\vskip4pt  

This completes the proof of Theorem~\ref{thm:recursive}.
\Endproof\vskip4pt

We now demonstrate that the UTP is well-behaved under connected sums.

\demo{\scshape Theorem~\ref{thm:connectedsum}} 
{\it If two knot types $\K_1$ and $\K_2$ satisfy the {\rm UTP,} then their
connected sum $\K_1\# \K_2$ satisfies the} UTP.

\Proof 
The following is the key claim:

\begin{claim}
Every solid torus $N$ with convex boundary which represents
$\K_1\#\K_2$ can be thickened to a standard neighborhood $N'$ of a
Legendrian curve in $\mathcal{L}(\K_1\#\K_2)$. 
\end{claim} 

\Proof 
Applying the Giroux Flexibility Theorem, $\bdry N$ can be put in {\em
standard form}, with meridional Legendrian rulings.  Let $S$ be the
separating sphere for $\K_1\#\K_2$ --- we arrange $S$ so it (1) is
convex, (2) intersects $N$ along two disks, and (3) intersects $\bdry
N$ in a union of Legendrian rulings.  Moreover, on the annular portion
of $S\setminus (\K_1\#\K_2)$, we may assume that (4) there are no
$\bdry$-parallel arcs, since otherwise $N$ can be thickened further by
attaching the corresponding bypasses. Now, cutting $S^3$ along $S$ and
gluing in copies of the standard contact $3$-ball $B^3$ with convex
boundary, we obtain solid tori $N_i$, $i=1,2$, (with convex boundary)
which represent $\K_i$.   

Since $\K_i$ satisfies the UTP, there exists a thickening of $N_i$ to
$N_i'$, where $N_i'$ is the standard neighborhood of a Legendrian knot
$L_i\in \L(\K_i)$.  Also arranging $\bdry N_i'$ so that it admits
meridional Legendrian rulings, we take an annulus from a Legendrian
ruling $\gamma'_i$ on $\bdry N_i'$ to a Legendrian ruling $\gamma_i$
on $\bdry N_i\cap \bdry N$.  If $\tb (\gamma_i)<-1$, then the
Imbalance Principle, together with the fact that $\tb(\gamma_i')=-1$,
yields enough bypasses which can be attached onto $\bdry N_i$ to
thicken $N_i$ into the standard neighborhood of a Legendrian knot.   

However, upon closer inspection, it is evident that the bypasses
produced can be attached onto $N$ inside the original $S^3$.  This
produces a thickening of $N$ to $N'$, which has boundary slope
${1\over m}$ (i.e., is the standard neighborhood of a Legendrian knot
in $\L(\K)$) measured with respect to $\mathcal{C}_{\K_1\#\K_2}$.   
\Endproof\vskip4pt  

Condition~1 of the UTP follows immediately from the claim.  To prove
Condition~2, we need to show that a standard 
neighborhood $N'$ of a Legendrian knot in $\L(\K_1\# \K_2)$ can be
thickened to $N''$ which is the standard neighborhood of a maximal
$\tb$ representative of $\L(\K_1\#\K_2)$.  This is equivalent to
showing any Legendrian knot $L'$ in $\L(\K_1\#\K_2)$ can be
destabilized to a maximal $tb$ representative. Given $L'\in
\L(\K_1\#\K_2)$, then $L'$ can be written as $L'_1\#L'_2$, with
$L_i'\in \L(\K_i)$, $i=1,2$.  Each $L_i'$ can be destabilized to a
maximal $\tb$ representative $L_i''$ by the UTP for each $\K_i$.
Since  
$$\overline{\tb}(\K_1\#\K_2) = \overline{\tb}(\K_1) +
\overline{\tb}(\K_2) +1,$$ 
by \cite{EH2}, we simply take $L''= L_1'' \# L_2''$.  This proves 
Theorem~\ref{thm:connectedsum}.
\hfill\qed

%%%%%%%%%%%%%%%%%%%%%%%%%%%%%%%%%%%%%%%%%%%%%%%%%%%%%%%
\section{Non-uniformly-thick knots and non-destabilizability} 
\label{section3}
%%%%%%%%%%%%%%%%%%%%%%%%%%%%%%%%%%%%%%%%%%%%%%%%%%%%%%%

We prove the following more precise version of 
Theorem~\ref{thm: positive torus}.

\demo{\scshape Theorem~\ref{thm: positive torus}} {\it The $(2,3)$-torus knot 
does not satisfy the} UTP.
\Enddemo
 
Although our considerations will work for any $(p,q)$-torus knot with 
$q>p>0$,  we assume for simplicity that $\K$ is a $(2,3)$-torus
knot, in order to keep the arguments simpler in a few places.

\Proof 
The goal is to exhibit solid tori $N$ representing $\K$, which cannot
be thickened to the maximal thickness.  The overall strategy is not
much different from the strategy used in \cite{EH3} and \cite{EH4} to
classify and analyze tight contact structures on Seifert fibered
spaces over $S^2$ with three singular fibers.  The plan is as follows:
we work backwards by starting with an arbitrary solid torus $N$ which
represents $\K$ and attempting to thicken it.  This gives us a list
$N_k$ of potential non-thickenable candidates, as well as tight
contact structures on their complements $S^3\setminus N_k$
(Lemma~\ref{slopes}).  In Lemma~\ref{tight verification} we prove that
the decomposition into $N_k$ and $S^3\setminus N_k$ actually exists
inside the standard tight $(S^3,\xi_{std})$ and in Lemma~\ref{no
  thickening} we prove the $N_k$ indeed resist thickening.  

Let $T$ be an oriented standardly embedded torus in $S^3$ which bounds 
solid tori $V_1$ and $V_2$ on opposite sides and which contains a 
$(2,3)$-torus knot $\K$.   Suppose $T=\bdry V_1$ and $T=-\bdry 
V_2$.  Also let $F_i$, $i=1,2$, be the core curve for $V_i$.  In 
\cite{EH1} it was shown that $\overline{\tb}(\K)=pq-p-q= 1$.  
Measured with respect to the coordinate system $\mathcal{C}'_{F_i}$
for either $i$, $\overline{t}(\K,\mathcal{C}'_{F_i})=-p-q=-5$, which 
corresponds to a slope of $-{1\over 5}$.

\begin{lemma}\label{slopes}
Suppose the solid torus $N$ representing $\K$ resists thickening.
Then ${\rm slope} (\Gamma_{\bdry N})=-{k+1\over 6k+5}${\rm ,} where $k$ is
a nonpositive integer and the slope is measured with respect to
$\mathcal{C}'_{F_i}$. 
\end{lemma}

\Proof 
Let $L_i$, $i=1,2$, be a Legendrian representative of $F_i$ with 
Thurston-Bennequin invariant $-m_i$, where $m_i>0$. If $N(L_i)$ is the
standard neighborhood of $L_i$, then ${\rm slope} (\Gamma_{\bdry
  N(L_i)})=-{1\over m_i}$ with respect to the coordinate system
$\mathcal{C}_{F_i}$.  We recast these slopes with respect to a new  
coordinate system $\mathcal{C}$ which identifies 
$T\stackrel\sim\rightarrow \R^2/\Z^2$, where $\K$ (viewed as sitting 
on $T$) corresponds to $(0,1)$. 

First we change coordinates from $\mathcal{C}_{F_1}$ to $\mathcal{C}$.
Consider the oriented basis $((2,3), (1,2))$ with respect to
$\mathcal{C}_{F_1}$; we map it to $((0,1),(-1,0))$ with respect to
$\mathcal{C}$.  This corresponds to the map $A_1=\begin{pmatrix} 3 &
-2 \\ 2 & -1 \end{pmatrix}.$   (Here we are viewing the vectors as
column vectors and multiplying by $A_1$ on the left.)  Then $A_1$ maps
$(-m_1,1)\mapsto (-3m_1-2, -2m_1-1)$.  Since we are only interested in
slopes, let us write it instead as $(3m_1+2, 2m_1+1)$.      

Similarly, we change from $\mathcal{C}_{F_2}$ to $\mathcal{C}$.  The
only thing we need to know here is that $(-m_2,1)$ with respect to
$\mathcal{C}_{F_2}$ maps to $(2m_2+3, m_2+2)$ with respect to~$\mathcal{C}$. 

Given a solid torus $N$ which resists thickening, let $L_i$, $i=1,2$, 
be a Legendrian representative of $F_i$ which maximizes $\tb(L_i)$ in 
the complement of $N$ (subject to the condition that $L_1\sqcup L_2$ 
is isotopic to $F_1\sqcup F_2$ in the complement of $N$).  View 
$S^3\setminus (N(L_1)\cup N(L_2)\cup N)$ as a Seifert fibered space 
over the thrice-punctured sphere, where the annuli which connect among 
$N(L_1)$, $N(L_2)$, and $N$ admit fibrations by the Seifert fibers.  
Now suppose $3m_1+2\not = 2m_2+3$.  Then we apply the Imbalance 
Principle to a convex annulus $A'$ between $N(L_1)$ and $N(L_2)$ to 
find a bypass along $N(L_i)$.  This bypass in turn gives rise to a 
thickening of $N(L_i)$, allowing the increase of $\tb(L_i)$ by one.  
Eventually we arrive at $3m_1+2=2m_2+3$ and a convex annulus $A'$ 
which has no $\bdry$-parallel arcs (hence we may assume $A'$ is in {\em
  standard form}).  Moreover, the denominator 
of ${\rm slope} (\Gamma_{\bdry N})$ must also equal $3m_1+2=2m_2+3$,
since otherwise $N$ admits a thickening.  
Since $m_i>0$, the smallest solution to $3m_1+2=2m_2+3$ is 
$m_1=1$, $m_2=1$.  All the other positive integer solutions are 
therefore obtained by taking $m_1=2k+1$, $m_2=3k+1$, with $k$ a 
nonnegative integer.  

We now compute the slope of the dividing curves on $\bdry ( N(L_1)\cup 
N(L_2)\cup N(A'))$, measured with respect to 
$\mathcal{C}_{F_1}'=\mathcal{C}_{F_2}'$, after edge-rounding.  Here 
$N(A')$ stands for the $I$-invariant neighborhood of the convex 
annulus $A'$.  We have: $$ -{2m_1+1\over 3m_1+2} + {m_2+2\over 2m_2+3} 
-{1\over 6k+5}= -{ 4k+3\over 6k+5} + {3k+3\over 6k+5} -{1\over 6k+5}=- 
{k+1\over 6k+5}.$$ For small $k$ we get $-{1\over 5}< -{2\over 11}< 
-{3\over 17}< -{4\over 23}<\dots < -{1\over 6}$.
\Endproof\vskip4pt  

Let $N_k$ be a tight solid torus representing $\K$ so that the
boundary slope is $-{k+1\over 6k+5}$ with respect to
$\mathcal{C}_{F_i}'$ and $\#\Gamma_{N_k}=2$.  (There are exactly two
tight contact structures on $N_k$ which satisfy the given boundary
conditions, and they are both universally tight.) Let
$M_k=S^3\setminus N_k.$  From the above discussion, if $N_k$ is to
resist thickening, then we know that $M_k$ must be contactomorphic to
the manifold obtained from $N(L_1)\cup N(L_2)$ by adding a standard
neighborhood of a convex annulus $A'.$  $M_k$ is a
Seifert fibered space and has a degree 6 cover $\widetilde{M_k}$
diffeomorphic to $S^1$ times a punctured torus (cf.\ \cite{EH1}). One may
easily\break\vskip-12pt\noindent check that the pullback of the tight contact structure to
$\widetilde{M_k}$ admits an isotopy where the $S^1$ fibers become
Legendrian and have twisting number $-(6k+5)$ with respect to the
product framing.  

\begin{lemma} \label{tight verification}
The standard tight contact structure on $S^3$ splits into a\break
\/{\rm (}\/universally\/{\rm )}\/ tight contact structure on $N_k$ with boundary slope
$-{k+1\over 6k+5}$ and the tight contact structure on $M_k$ described
above.  
\end{lemma}

\Proof 
Let $N_k$ be a (universally) tight solid torus described above and let
$A$ be a convex annulus in standard form from $N_k$ to itself, such
that the complement of $R=N_k\cup N(A)$ in $S^3$ consist of standard
neighborhoods $N(L_i)$, $i=1,2$.  Here $N(A)$ is the $I$-invariant neighborhood
of $A$.  (Observe that $R$ is also contact isotopic to $N_k \cup N(A')$.) 
 
For either choice of contact structure on $N_k$, the contact structure
on $R$ can be isotoped to be transverse to the fibers of $R$ (where
the fibers are parallel to $\K$), while preserving the dividing set on
$\bdry R$.  Such a {\em horizontal} contact
structure is universally tight.  (For more details of this standard
argument, see for example \cite{H2}.)
 
Once we know that the contact structure on $R$ is tight, we just need
to apply the classification of tight contact structures on solid tori
and thickened tori.  In fact, any tight contact structure on
$R=T^2\times [1,2]$ with boundary conditions $\#\Gamma_{T_1}=\#\Gamma_{T_2}=2$ and
${\rm slope} (\Gamma_{T_1})=-{1\over m_1}$,
${\rm slope} (\Gamma_{T_2})=-m_2$ (here $m_i$ are positive integers)
glues together with 
$N(L_1)$ and  $N(L_2)$ to give the tight contact structure on $S^3$.  
\hfill\qed

\begin{lemma} \label{no thickening}
The tight solid torus $N_k$ does not admit a thickening to a solid 
torus $N_{k'}$ whose boundary slope is $-{k'+1\over 6k'+5}${\rm ,} where 
$k'< k$.   More generally{\rm ,} $N_k$ does not admit any nontrivial 
thickenings{\rm ,} i.e.{\rm ,} no thickenings with a boundary slope different from 
$-{k+1\over 6k+5}$. 
\end{lemma}   

\Proof 
If $N_k$ can be thickened to $N_{k'}$, then there exists a Legendrian
curve isotopic to the regular fiber of the Seifert fibered space
$M_k=S^3\setminus N_k$ with twisting number $> -(6k+5)$, measured with
respect to the Seifert fibration.  (Take a ruling curve on $\partial
N_{k'}\subset M_k.$) Pulling back to the sixfold cover
$\widetilde{M_k}$, we have a Legendrian knot which is topologially
isotopic to a fiber but has twisting number $> -(6k+5).$  However, we
claim that the maximal twisting number for a fiber in
$\widetilde{M_k}$ is $-(6k+5).$  One way to see this is to add a solid
torus to $\widetilde{M_k}$ to obtain $T^3$ and extend the 
contact structure so that all the $S^1$ fibers in $T^3$ are Legendrian
with twisting $-(6k+5)$.  We can now apply the classification of tight
contact structures on $T^3$ due to Giroux and Kanda (see \cite{K}) to
conclude that the maximal twisting number for a fiber is $-(6k+5).$

Next, suppose $N_k$ admits a nontrivial thickening $N'$ (not
necessarily of type $N_{k'}$).  Then we use 
the argument in Lemma~\ref{slopes} to find 
Legendrian curves $L_i\subset S^3\setminus N'$ which maximize the 
twisting number amongst Legendrian curves isotopic to $F_i$ in 
$S^3\setminus N'$, and a convex annulus from $N(L_1)$ to $N(L_2)$, so 
that $\bdry (N(L_1)\cup N(L_2)\cup N(A'))$ has some slope $-{k'+1\over 
6k'+5}$, $k'<k$.  This puts us in the case treated in the previous 
paragraph.   
\Endproof  

This completes the proof of Theorem~\ref{thm: positive torus}. 
\Endproof 

As a corollary of the above investigation we have: 

\demo{\scshape Theorem~\ref{thm: no destabilization}} {\it Let $\K'$ be the 
$(2,3)$-cable of the $(2,3)$-torus knot $\K$.  Then there exists a 
Legendrian knot $L\in\mathcal{L}(\K')$ which does not admit any 
destabilization, yet satisfies $\tb(L) < \overline{\tb}(\K')$.}

\Proof    
Let $N_k$ be a solid torus which resists thickening; say $k=1$.  
Then the boundary slope of $N_1$ is $-{2\over 11}$, measured with 
respect to $\mathcal{C}'_{F_i}$.  We choose a slope $-{a\over 
b} <-{2\over 11}$ whose corresponding simple closed 
curve, denoted $(-b,a)$, has fewer intersections with the simple 
closed curve $(-11,2)$ than with any other simple closed curve whose 
corresponding slope $-{c\over d}$ satisfies $-{2\over 11} < -{c\over 
d} < 0$.  To verify that $-{a\over b}=-{3\over 16}$ works, consider 
the standard Farey tessellation of the hyperbolic unit disk.  
Since there mutually are edges among $-{1\over 5}< -{3\over 16} 
<-{2\over 11}$,\break\vskip-11pt\noindent $-{a\over b}=-{3\over 16}$ is shielded from any 
$-{c\over d} > -{2\over 11}$ by the edge \pagebreak from $-{1\over 5}$ to 
$-{2\over 11}$.  Therefore, to get from $-{3\over 16}$ to $-{c\over 
d}$ we need at least two steps, implying that $(-16,3)$ and $(-11,2)$ 
have fewer  intersections than $(-16,3)$ and $(-d,c)$.  Now, by 
changing coordinates from $\mathcal{C}'_{F_i}$ to $\mathcal{C}_{\K}$, 
we see that the slope $-{a\over b}=-{3\over 16}$ corresponds to the 
$(2,3)$-cable of the $(2,3)$-torus knot.  

First observe that there is a Legendrian knot $L'\in \L(\K')$ which
sits inside the solid torus $N_0$ with ${\rm slope} (\Gamma_{\bdry
  N_0})=-{1\over 5}$ (with respect to $\mathcal{C}'_{F_i}$), as a
Legendrian divide on a convex torus which 
is isotopic to (but not contact isotopic to) $\bdry N_0$ and which has
slope $-{3\over 16}$.  By the classification of tight contact
structures on solid tori, such a convex torus exists because $-{3\over
  16}> -{1\over 5}$.  This proves that
$\overline{t}(\K',\mathcal{C}'_\K)\geq 0$.  

Next we exhibit $L\in \L(\K')$ which cannot be destabilized to twisting  
number $0$ with respect to $\mathcal{C}'_{F_i}$.  Let $L$ be a
Legendrian ruling curve on $\bdry N_1$,  
where the ruling is of slope $-{3\over 16}$.  By 
construction, the twisting number $t(L,\mathcal{C}'_\K)=-1$, computed 
by intersecting $(-11,2)$ and $(-16,3)$.  

\begin{lemma}
$L$ cannot be destabilized.
\end{lemma}

\Proof 
The proof is an application of the state transition technique 
\cite{H3}.  Suppose that $L$ admits a destabilization.  Then there 
exists a convex torus $\Sigma$ isotopic to $\bdry N_1$ which contains 
$L$ as well as a bypass to $L$.  More conveniently, instead of 
isotoping both $L$ and the torus, we fix $L$ and isotop the torus from 
$\bdry N_1$ to $\Sigma$.  Then the annulus $B_0=(\bdry N_1)\setminus L$  
is isotoped to $B=\Sigma\setminus L$ relative to the boundary.  
Observe that $\Gamma_{B_0}$ consists of two parallel nonseparating 
arcs.  To get to $B$, we perform isotopy discretization, i.e., a 
sequence of bypass moves (which may well be trivial bypass 
attachments).  There can be no nontrivial bypasses attached onto $B_0$ 
from the exterior of $N_1$, since $N_1$ has maximal thickness.  

We claim there are no nontrivial bypasses from the interior as well.
First of all, since there are no Legendrian knots isotopic to $L$ with
twisting number zero inside $N_1$, no $\bdry$-parallel dividing curves
(on $B$) can be created by attaching bypasses from the interior.  On
the other hand, the slope (or holonomy) of the two separating arcs on
$B_0$ cannot be changed, since the only slope $-{c\over d}$ with
$-{c\over d} \geq -{2\over 11}$ with an edge (in the Farey
tessellation) to $-{3\over 16}$ is $-{2\over 11}$.  This proves that
all the state transitions for $B_0$ are trivial state transitions.  We
are unable to reach $B$.  
\Endproof\vskip4pt                
This completes the proof of Theorem~\ref{thm: no destabilization}.
\Endproof\vskip4pt

\section{Non-transverse-simplicity}\label{sec:nts}

\demo{\scshape Theorem~\ref{thm: transverse}} {\it Let $\K'$ be the 
$(2,3)$-cable of the $(2,3)$-torus knot $\K$.  Then $\K'$ is not 
transversely simple.}
\Enddemo

We first gather some preliminary lemmas.

\begin{lemma}   \label{tbcalc}
$\overline{\tb}(\K')=w(\K')=6$.
\end{lemma}

The proof of this lemma is identical to that of
Theorem~\ref{thm:torusknots}.  

\begin{lemma}   \label{max}
There are precisely two maximal Thurston-Bennequin representatives in
$\L(\K')${\rm ,} which we call $K_\pm$ and which have $tb(K_\pm)=6${\rm ,}
$r(K_\pm)=\pm 1$.  
\end{lemma}

\Proof 
Any $K\in \L(\K')$ with $\tb(K)=6$ can be realized as a Legendrian
divide on the boundary of a solid torus $N$ representing $\K$.  
By Lemma~\ref{slopes}, $N$ can be thickened to a solid torus $N'$
with ${\rm slope} (\Gamma_{\bdry N'})= -{1\over 5}$, measured with
respect to $\mathcal{C}'_{F_i}$.  This means that there are two  
possible tight contact structures on $N$, both universally tight, and 
the extension to $N'$ is determined by the tight contact structure on 
$N$.  Once $N'$ is determined, the tight contact structure on 
$S^3\setminus N'$ is unique up to isotopy, since $N'$ is the standard 
tubular neighborhood of the unique maximal $\tb$ representative of 
$\K$.   This proves that there are at most two maximal $\tb$ 
representatives of $\L(\K)$.

We now show that there are indeed two representatives by computing
their rotation numbers to be $r(K)=\pm 1$ (and hence showing they are
distinct).  To use Lemma~\ref{rcomp}, we need to know the rotation
number of a ruling curve $\lambda$ isotopic to $\K$ on $\partial N$
and the rotation number of a meridional ruling curve $\mu$ on
$\partial N.$  A ruling curve isotopic to $\K$ on $\partial N'$ has
rotation number 0 (by the Bennequin inequality). The region $R$
between $\partial N$ and $\partial N'$ (in $\CC_\K$ coordinates) has
relative half-Euler class 
$$PD(e(\xi),R)= \pm ((1,1)-(2,3))= \pm (-1,-2).$$ 
So $r(\lambda)=\mp 1.$  One similarly sees that $r(\mu)=\pm 2.$
Thus 
\vglue12pt
\hfill $r(K)= 2(\pm 2) + 3(\mp 1)= \pm 1. $
\Endproof\vskip4pt  

\begin{lemma}\label{L pm}
The only non-destabilizable representatives of $\L(\K')$ besides 
those which attain $\overline{\tb}(\K')$ are $L_\pm$ which have 
$\tb(L_\pm)=5$ and $r(L_\pm)=\pm 2$.  They are realized as Legendrian 
ruling curves on a convex torus isotopic to $T$ with dividing curves 
of slope $-{2\over 11}$ \/{\rm (}\/with respect to $\mathcal{C}'_{F_i}${\rm ),} and
which does not admit a thickening.   
\end{lemma}

\Proof 
Let $K$ be a non-destabilizable representative of $\L(\K')$.  Since 
$\overline{\tb}(\K')=6$, we can always place $K$ on the (convex) boundary
$\Sigma=\bdry N$ of a solid torus $N$ representing $\K$.  If $K$ is a
Legendrian divide on $\Sigma$, then we are in the case of
Lemma~\ref{max}.  If $K$ is not a Legendrian divide, then $K$ must intersect
$\Gamma_\Sigma$ efficiently, and we may assume that  $K$ is a
Legendrian ruling curve on $\Sigma$.  Slopes of $\Sigma$ will usually
be measured with respect to $\mathcal{C}'_{F_i}$.

We now show that if $s={\rm slope} (\Gamma_\Sigma)\not=-\frac{2}{11}$,
then $K$ can be destabilized (contradicting our assumption).  Note
that $s$ must be in $[-\frac{1}{5},0)$ and $s= -\frac{3}{16}$
corresponds to the situation in Lemma~\ref{max}.  In the following
cases, we find a convex torus $\Sigma'$ isotopic to and disjoint from
$\Sigma$ so that a simple closed curve of slope 
$-{3\over 16}$ has smaller geometric intersection with
$\Gamma_{\Sigma'}$ than with $\Gamma_\Sigma$.  The
destabilization is then a consequence of the Imbalance Principle.  If
$s\in[-{1\over 5},-{3\over 16})$, then there is
$\Sigma'\subset N$ with ${\rm slope} (\Gamma_{\Sigma'})=-{3\over 16}$. 
If $s\in(-{3\over 16}, -{2\over 11})$, then there exists
$\Sigma'$ of slope $-{3\over 16}$ outside $N$ (since $N$ can be
thickened to maximal width by Lemma~\ref{slopes}).  Similarly, if $s\in
(-{1\over 6},0)$, then there exists a $\Sigma'$ with  
${\rm slope} (\Gamma_{\Sigma'})=-{1\over 6}$, by using Lemma~\ref{slopes}. 
Next, if $s\in (-{2\over 11},-{1\over 6})$, there exists a $\Sigma'$
of slope $-{1\over 6}$ inside $N$ (it is not difficult to see that
this $\Sigma'$ works by referring to the Farey tessellation).
Therefore we are left with $s =-{1\over 6}$.   But then we use the 
classification of $\L(\K)$ to deduce that $N$ can be thickened to $N'$ 
with boundary slope $-{1\over 5}$, corresponding to a representative 
of $\L(\K)$ of maximal Thurston-Bennequin invariant.  We can now 
compare $\Sigma$ with $\Sigma'$ of slope $-{3\over 16}$ to 
destabilize. This proves that the only two places where we get 
stuck and cannot destabilize are $-{2\over 11}$ and $-{3\over 16}$. 

Now let $L\in \L(\K')$ be non-destabilizable representatives with 
$\tb(L)=5$.  Then they are Legendrian ruling curves on the boundary of 
a solid torus $N_1$, where ${\rm slope} (\Gamma_{\bdry N_1})=-{2\over 
11}$ with respect to $\mathcal{C}'_{F_i}$.  There are two possible 
tight contact structures on $N_1$, and they are both universally 
tight.  Since the tight contact structures on their complements 
$S^2\setminus N_1$ are always contact isotopic, there are at most two 
non-destabilizable, nonmaximal representatives.  Using 
Lemma~\ref{rcomp}, we obtain:
$$r(L) = 2(\pm 1) + 3(0)=\pm 2.$$
(Since $\mu$ intersects $\Gamma_{\bdry N_1}$ in four points,
$r(\mu)=\pm 1.$ It is also not hard to compute $r(\lambda)=0$ by using
the fact that there are no $\bdry$-parallel arcs on the Seifert
surface for $\lambda$.) 
Therefore $L_+$ and $L_-$ are distinguished by the contact 
structures on the solid torus $N_1$.
\hfill\qed

\begin{lemma}
$S_-(L_-)= S_-^2(K_-)$ and $S_+(L_+)=S_+^2 (K_+)$.
\end{lemma}

\Proof 
Since $L_-$ is a Legendrian ruling curve on $N_1$ with 
${\rm slope} (\Gamma_{\bdry N_1})=-{2\over 11}$, $S_-(L_-)$ is a 
Legendrian ruling curve on $\partial N_1'\subset N_1$, where $N'_1$ is a 
solid torus representing $\K$ and ${\rm slope} (\Gamma_{\bdry N'_1})= -{1\over 6}$.  
Similarly, since $K_-$ is a Legendrian divide on $N$ with
${\rm slope} (\Gamma_{\partial N})=-\frac{3}{16}$, $S^2_-(K_-)$
is a Legendrian ruling curve on $\partial N'\subset N$, where $N'$ is a 
solid torus representing $\K$ and ${\rm slope} (\Gamma_{\bdry N'})=-{1\over 6}$.
Now $N'$ and $N'_1$ are neighborhoods of Legendrian knots in
$\L(\K)$ with $tb=0.$  If the associated rotation numbers are
the same, then they are contact isotopic (by the Legendrian simplicity
of the $(2,3)$-torus knot).  One may easily check that 
the rotation numbers are indeed the same.  Therefore, there is an ambient
contact isotopy taking $N'$ to $N_1'$, and it simply remains to Legendrian 
isotop $S_-(L_-)$ to $ S_-^2(K_-)$ through ruling curves.
\Endproof\vskip4pt  

We are now ready to proceed with the proof of Theorem~\ref{thm: 
transverse}. 

\demo{Proof of Theorem~{\rm \ref{thm: transverse}}} 
In view of Theorem~\ref{stably}, it suffices to show that $S_+^k 
(L_-)$ is never equal to $S_+^k S_-(K_-)$ for all positive integers 
$k$ (and likewise $S_-^k (L_+)$ is never equal to $S_-^k S_+(K_+)$).  

Throughout this proof we use coordinates $\mathcal{C}'_{\K}$, unless
otherwise stated.  As above, let $N_1$ be a solid torus which
represents $\K$, does not admit a thickening, and has boundary
$\Sigma_0=\bdry N_1$, where $\#\Gamma_{\Sigma_0}=2$ and ${\rm slope} 
(\Gamma_{\Sigma_0}) = -\frac{2}{11}$.  Assuming we have already chosen
the correct $N_1$ (there were two choices), place the knot
$L=S_+^k(L_-)$ on $\Sigma_0$ as follows: 
if $A_0=\Sigma_0\setminus L$, then there are $k$ negative 
$\bdry$-parallel arcs on the left-hand edge $L_l$ of $A_0$ and $k$ 
positive $\bdry$-parallel arcs on the right-hand edge $L_r$ of $A_0$.  
Here $A_0$ is oriented so that $\bdry A_0= L_r - L_l$, where $L_r$ and 
$L_l$ are oriented copies of $L$.  (The sign of a $\bdry$-parallel 
arc is the sign of the region it cuts off.) See Figure~\ref{initial}
for a possible $\Gamma_{A_0}$.  When we draw annuli, we 
will usually present rectangles, with the understanding that the
top and the bottom are identified. \setcounter{figure}{2}
\begin{figure}[ht]
{\epsfysize=2in \centerline{\epsfbox{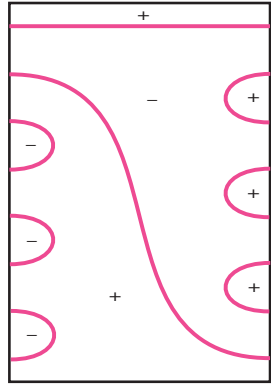}}
\caption{The ``initial configuration'' $\Gamma_{A_0}$.  The left-hand boundary 
is $L_l$ and the right-hand boundary is $L_r$.  They glue to give 
$\Sigma_0$.}    
\label{initial}} 
\end{figure}
 
%For convenience when drawing pictures we will use coordinates $\CC$ on $T$ 
%where $\pm(1,0)$ and $\pm(0,1)$ correspond to slopes 
%${-\frac{2}{11}}$ and $-{3\over 16}$ with respect to
%$\mathcal{C}'_\K$.

The key claim is the following:

\begin{claim}\label{same slope}
Every convex torus which contains $L$ and is isotopic to $\Sigma_0$
has slope $-\frac{2}{11}$.
\end{claim}

This would immediately show that $S_+^k(L_-)$ is never equal to $S_+^k 
S_-(K_-)$. To prove this fact, we use the state traversal technique.
If $\Sigma$ also contains $L$ and is isotopic to $\Sigma_0$ (not
necessarily relative to $L$), then we can use the standard properties
of incompressible surfaces in Seifert fibered spaces to conclude that
$\Sigma$ must be isotopic to $\Sigma_0$ {\em relative to $L$}.  Therefore,
it suffices to show that the slope of the dividing set does not change 
under any isotopy of $\Sigma_0$ relative to $L$.  Although we would
like to say that the isotopy leaves the dividing set of $\Sigma_0$
invariant, this is not quite true.  It is not difficult to see (see
Figure~\ref{inductive}) that the number of dividing curves can increase,
although the slope should always remain the same according to
Claim~\ref{same slope}.  Starting with
$\Sigma=\Sigma_0$, we inductively assume the following:

\begin{figure}[ht]
{\epsfysize=2in \centerline{\epsfbox{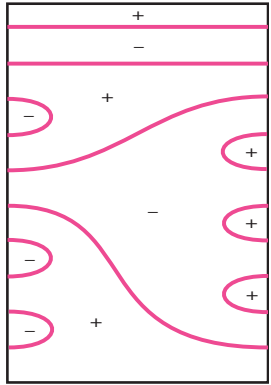}}
\caption{A potential $\Sigma$ in the inductive step.}   
\label{inductive}} 
\end{figure}

\demo{Inductive hypothesis} 
\begin{enumerate}
\item $\Sigma$ is a convex torus which contains $L$ and satisfies
  $2\leq \#\Gamma_\Sigma \leq 2k+2$ and ${\rm slope}  (\Gamma_\Sigma)
  = -\frac{2}{11}$.  
\item $\Sigma$ is ``sandwiched'' in a $[0,1]$-invariant
   $T^2\times[0,1]$ with ${\rm slope}  (\Gamma_{T_0}) = {\rm slope} 
   (\Gamma_{T_1}) = 
   -\frac{2}{11}$ and $\#\Gamma_{T_0} = \#\Gamma_{T_1}=2$.  (More
   precisely, $\Sigma\subset T^2\times (0,1)$ and is parallel to
   $T^2\times \{i\}$.)
\item There is a contact diffeomorphism
  $\phi:S^3\stackrel\sim\rightarrow S^3$ which takes $T^2\times [0,1]$
  to a standard $I$-invariant neighborhood of $\Sigma_0$ and matches
  up their complements. 
\end{enumerate}

Suppose we isotop $\Sigma$ relative to $L$ into another convex torus
$\Sigma'$.  Then the standard state traversal machinery \cite{H3}
implies that we may assume that the isotopy is performed in discrete
steps, where each step is given by the attachment of a bypass.
$\Sigma$ bounds a solid torus $N$ on one side, and we say that the
bypass is attached ``from the inside'' or ``from the back'' if the
bypass is in the interior of $N$ and the bypass is attached ``from the
outside'' or ``from the front''  if the bypass is in the exterior of
$N$. (Also for convenience assume that $T_0$ is inside $N$ and $T_1$
is outside $N$.)  We prove the inductive hypothesis still holds after
all existing bypass attachments.  

\begin{lemma}\label{nonbhd}
The Legendrian knot $L$ cannot sit on a convex torus $\Sigma$ in $N_1$
that is isotopic to $\partial N_1$ and satisfies $\#\Gamma_{\Sigma}=2$ and
${\rm slope} (\Gamma_{\Sigma})=-\frac{1}{6}$.  
\end{lemma}

\Proof 
The convex torus $\Sigma$ bounds the standard neighborhood of a
Legendrian knot in $\L(\K)$ with $\tb=0$ and $r=-1$ (i.e.\, $S_-$ of the
maximal $\tb$ representative of $\L(\K)$).  Computing as in
Lemma~\ref{L pm}, we find that a Legendrian ruling curve of slope
$-{3\over 16}$ on $\Sigma$ must be $S_-(L_-)$.  Therefore, if
$L\subset \Sigma$, then $L$ must be a stabilization of
$S_-(L_-)$. However, this contradicts the fact that $L= S_+^k(L_-)$ by
a simple $(r,\tb)$-count. 
\Endproof\vskip4pt  

\begin{lemma}
Given a torus $\Sigma$ satisfying the inductive hypothesis{\rm ,} any bypass
attached to $A=\Sigma\setminus L$ will not change the slope of the
dividing set. 
\end{lemma}

\Proof 
If the bypass is attached from the outside, then the slope cannot
change or this would give a thickening of our non-thickenable solid
torus.  If the bypass is attached from the inside, then let $\Sigma'$ be
the torus obtained after the bypass is attached.  By examining the
Farey tessellation, we see that $s={\rm slope} (\Gamma_{\Sigma'})$
must lie in $[-\frac{2}{11}, -\frac{1}{6}]$. Since Lemma~\ref{nonbhd}
disallows $s=-{1\over 6}$, suppose that $s\in (-{2\over 11}, -{1\over
  6})$.  Let $\Sigma''$ be a convex torus of slope $-{1\over 6}$ and
$\#\Gamma=2$ in the interior of the solid torus bounded by $\Sigma'$.
Take a Legendrian curve $L'$ on $\Sigma'$ which is parallel to and
disjoint from $L$, and intersects $\Gamma_{\Sigma'}$ minimally.
Similarly, consider $L''$ on $\Sigma''$.  Using the Farey
tessellation, it is clear that $|\Gamma_{\Sigma'}\cap L'| >
|\Gamma_{\Sigma''}\cap L''|$.  Thus the Imbalance Principle gives bypasses
for $\Sigma'$ that are disjoint from $L$.  After successive attachments
of such bypasses, we eventually obtain $\Sigma'''$ of slope $-{1\over
  6}$ containing $L$, contradicting Lemma~\ref{nonbhd}.  
\Endproof\vskip4pt  

Therefore we see that Condition (1) is preserved.

\begin{lemma}
Given a torus $\Sigma$ satisfying the inductive hypothesis{\rm ,} any bypass
attached to $A$ will preserve Conditions {\rm (2)} and {\rm (3)}.
\end{lemma}

\Proof 
Suppose $\Sigma'$, is obtained from $\Sigma$ by a single bypass move.
We already know that ${\rm slope}  (\Gamma_{\Sigma'}) = {\rm slope} 
(\Gamma_{\Sigma})$, and, assuming the bypass move was not trivial,
$\#\Gamma$ is either increased or decreased by 2.  Suppose first that 
$\Sigma'\subset N$, where $N$ is the solid torus bounded by
$\Sigma$. For convenience, suppose $\Sigma=T_{0.5}$ inside
$T^2\times[0,1]$ satisfying Conditions (2) and (3) of the inductive
hypothesis. Then we form the new $T^2\times[0.5,1]$ by taking the old  
$T^2\times[0.5,1]$ and adjoining the thickened torus between $\Sigma$ 
and $\Sigma'$.  Now, $\Sigma'$ bounds a solid torus $N'$, and, by the 
classification of tight contact structures on solid tori, we can 
factor a nonrotative outer layer which is the new $T^2\times[0,0.5]$.  

On the other hand, suppose $\Sigma'\subset (S^3\setminus N)$.  We 
prove that there exists a nonrotative outer layer $T^2\times[0.5,1]$ 
for $S^3\setminus N'$, where $\#\Gamma_{T_1}=2$.  This follows from 
repeating the procedure in the proof of Theorem~\ref{thm: positive 
torus}, where Legendrian representatives of $F_1$ and $F_2$ were 
thickened and then connected by a vertical annulus  --- this time the 
same procedure is carried out with the provision that the 
representatives of $F_1$ and $F_2$ lie in $S^3\setminus N'$.  Once the 
maximal thickness for representatives of $F_1$ and $F_2$ is obtained, 
after rounding we get a convex torus in $S^3\setminus N'$ parallel to 
$\Sigma'$ but with $\#\Gamma=2$.  Therefore we obtain a nonrotative 
outer layer $T^2\times[0.5,1]$. 
\Endproof\vskip4pt  

This completes the proof of Theorem~\ref{thm: transverse}.
\hfill\qed

\demo{Acknowledgements}  
Part of this work was done during the Workshop on Contact Geometry in
the fall of 2000.   The authors gratefully acknowledge the support of
Stanford University and the American Institute of Mathematics during
this workshop.  The first author was supported by an NSF Postdoctoral Fellowship 
(DMS-9705949), NSF Grant  DMS-0203941, an NSF CAREER Award (DMS-0239600), and an
NSF FRG grant (FRG-0244663), and
the second author by an Alfred P.\ Sloan Research Fellowship, NSF
Grant  DMS-0072853, and an NSF CAREER Award (DMS-0237386). The second author also
thanks the University of Tokyo and the Tokyo Institute of Technology
for their hospitality during his stay in 2003, when the paper was completed.


\references {cables}

\bibitem[Be]{Be}
\name{D.\ Bennequin}, Entrelacements et \'equations de Pfaff,
Third Schnepfenried geometry conference, Vol.\ 1 (Schnepfenried, 1982),
87--161, {\it Ast{\hskip.5pt\rm \'{\hskip-5pt\it e}}risque\/} \textbf{107--108}, Soc.\ Math.\ France, Paris,
1983.

\bibitem[BM]{BM}
\name{J.\ Birman} and \name{W.\ Menasco}, Stabilization in the braid 
groups II: transversal simplicity of transverse knots, preprint 2002.

\bibitem[Ch]{Ch}
\name{Y.\ Chekanov}, Differential algebra of Legendrian links, 
{\it Invent.\ Math\/}.\ {\bf 150}  (2002), 441--483.

\bibitem[Co]{Colin}
\name{V.\ Colin}, Sur la stabilit\'e, l'existence et l'unicit\'e des
structures de contact en dimension $3$, Ph.\ D.\ Thesis, \'Ecole normale
sup\'erieure de Lyon, October 1998.

\bibitem[EGH]{EGH}
\name{Y.\ Eliashberg, A.\ Givental}, and \name{H.\ Hofer},   Introduction to 
symplectic field theory, GAFA 2000 (Tel Aviv,
1999), {\it Geom.\ Funct.\ Anal\/}.\ {\bf 2000}, 
Special Volume, Part II, 560--673. 

\bibitem[El]{Eliashberg92}
\name{Y.\ Eliashberg}, Contact $3$-manifolds twenty years since J.\ Martinet's 
work, {\it Ann.\ Inst.\ Fourier\/} ({\it Grenoble\/}) \textbf{42} 
(1992), 165--192.

\bibitem[EF]{EF}
\name{Y.\ Eliashberg} and \name{M.\ Fraser}, Classification of topologically trivial 
Legendrian knots, in \textit{Geometry}, \textit{Topology}, 
\textit{and Dynamics} (Montreal, PQ, 
1995), 17--51, {\it CRM Proc.\ Lecture Notes\/} \textbf{15}, A.\ M.\ S., 
Providence, RI, 1998.
 
\bibitem[EH1]{EH1}
\name{J.\ Etnyre} and \name{K.\ Honda}, Knots and contact geometry I: torus knots and 
the figure eight knot, {\it J.\ Symplectic Geom\/}.\ {\bf 1} (2001),  
63--120.

\bibitem[EH2]{EH2}
\name{J.\ Etnyre} and \name{K.\ Honda}, On connected sums and Legendrian 
knots, {\it Adv.\ Math\/}.\ {\bf 179} (2003), 59--74.

\bibitem[EH3]{EH3}
\bibline,  On the nonexistence of tight contact structures,  
{\it Ann.\ of Math\/}.\ {\bf 153} (2001), 749--766.

\bibitem[EH4]{EH4}
\bibline,  Tight contact structures with no symplectic fillings, 
{\it Invent.\ Math\/}.\ {\bf 148} (2002), 609--626.

\bibitem[Ga]{Gay}
\name{D.\ Gay},  Symplectic $2$-handles and transverse links, {\it Trans.\ 
Amer.\ Math.\ Soc\/}.\ {\bf 354} (2002), 1027--1047.

\bibitem[Gi1]{Gi1}
\name{E.\ Giroux}, Convexit\'e en topologie de contact, {\it Comment.\ Math.\ 
Helv\/}.\ \textbf{66} (1991), 637--677.

\bibitem[Gi2]{Gi2}
\bibline,  Structures de contact en dimension trois
et bifurcations des feuilletages de surfaces, {\it Invent.\ Math\/}.\ 
{\bf 141} (2000), 
615--689.

\bibitem[H1]{H1}
\name{K.\ Honda}, On the classification of tight contact
structures I, {\it Geom.\ Topol\/}.\ {\bf 4} (2000), 309--368.

\bibitem[H2]{H2}
\bibline,  On the classification of tight contact structures II,
{\it J.\ Differential Geom\/}.\ {\bf 55} (2000), 83--143.

\bibitem[H3]{H3}
\bibline, Gluing tight contact structures, {\it Duke Math.\ J\/}.\ {\bf 115} 
(2002), 435--478.

\bibitem[H4]{H4}
\bibline, Factoring nonrotative $T^2\times I$ layers, 
Erratum to ``On the classification of tight contact structures 
I'', {\it Geom.\ Topol\/}.\ {\bf 5} (2001), 925--938.


\bibitem[K]{K}
\name{Y.\ Kanda},  The classification of tight contact structures on
the $3$-torus, {\it Comm.\ Anal.\ Geom\/}.\ {\bf 5} (1997),  413--438.

\bibitem[M1]{Menasco01}
\name{W.~Menasco},  On iterated torus knots and transversal knots,  
{\it Geom.~Topol\/}.\ {\bf 5} (2001), 651--682.

\bibitem[M2]{Menasco?}
\bibline,  Erratum to: On iterated torus knots and transversal knots,
in preparation. 

\bibitem[Ng]{Ng}
\name{L.\ Ng}, Computable Legendrian invariants, {\it Topology\/} {\bf 42} 
(2002), 55--82.


\Endrefs
\end{document}